# THE LEFT LOCALIZATION PRINCIPLE, COMPLETIONS, AND COFREE $G$-SPECTRA

LUCA POL AND JORDAN WILLIAMSON

ABSTRACT. We show under mild hypotheses that a Quillen adjunction between stable model categories induces another Quillen adjunction between their left localizations, and we provide conditions under which the localized adjunction is a Quillen equivalence. Moreover, we show that in many cases the induced Quillen equivalence is symmetric monoidal. Using our results we construct a symmetric monoidal algebraic model for rational cofree $G$-spectra. In the process, we also show that $L$-complete modules provide an abelian model for derived complete modules.

## Contents

## 1. INTRODUCTION

In this paper we investigate the interplay between adjoint pairs and localizations. In homotopy theory there are two versions of localizations available: the left and right Bousfield localization. The former is ubiquitous in chromatic stable homotopy theory, while the latter has seen interesting applications in the study of torsion objects in algebraic categories, see [25, §5]. Often in the literature the right Bousfield localization is called cellularization since in the stable setting it picks out the localizing subcategory on the set of cells (if they are stable).

We now give an informal overview of our results and refer to the main body of the paper for the precise statements.

**The Cellularization Principle.** Let $\mathcal{C}$ be a stable model category, and let $\mathcal{K}$ be a set of objects of $\mathcal{C}$. The Cellularization Principle of Greenlees-Shipley [25] provides conditions under which a Quillen adjunction $F: \mathcal{C} \rightleftarrows \mathcal{D} : G$ descends to a Quillen equivalence

$$F : \operatorname{Cell}_{\mathcal{K}}\mathcal{C} \rightleftarrows \operatorname{Cell}_{F\mathcal{K}}\mathcal{D} : G$$

between the cellularizations. The Cellularization Principle is a crucial ingredient in the construction of algebraic models for rational equivariant spectra, see for instance [28]. There is also a version of the Principle where the cells are passed along the right adjoint, and a variant [9, §5.1] in which symmetric monoidal structures are taken into account. The main limitation of the Cellularization Principle is that the preservation of symmetric monoidal structures is *not* automatic.

Since the symmetric monoidal structure need not be preserved by cellularization, the symmetric monoidal version of the Cellularization Principle requires stronger assumptions. For instance, when passing cells along the right adjoint, the Cellularization Principle gives a symmetric monoidal Quillen equivalence between the cellularizations if the original adjunction was *already* a symmetric monoidal Quillen equivalence [9, 5.1.7].

On the other hand, the monoidal structure is often preserved by left Bousfield localization.

**The Left Localization Principle.** The Left Localization Principle which we develop, gives mild conditions under which a symmetric monoidal Quillen adjunction $F: \mathcal{C} \rightleftarrows \mathcal{D} : G$ descends to a *symmetric monoidal* Quillen equivalence between the homological localizations. For an object $E$ of a stable, symmetric monoidal model category $\mathcal{C}$, the homological localization $L_E \mathcal{C}$ is the localization of $\mathcal{C}$ at the class of $E$-equivalences, that is those morphisms that become equivalences after tensoring with $E$.

**Theorem** (3.13). *Let $\mathcal{C}$ and $\mathcal{D}$ be stable, symmetric monoidal model categories, $E$ an object of $\mathcal{C}$ and $F : \mathcal{C} \rightleftarrows \mathcal{D} : G$ be a symmetric monoidal Quillen adjunction. Suppose that $\mathcal{C}$ is homotopically compactly generated by a set $\mathcal{K}$ of objects and that $\mathcal{D}$ is homotopically compactly generated by $F\mathcal{K}$. Suppose that:*

(i) *The derived unit map $K \to GFK$ is an $E$-equivalence for all $K \in \mathcal{K}$;*
(ii) *$G$ sends $FE$-equivalences to $E$-equivalences.*

*Then the induced Quillen adjunction*

$$F : L_E \mathcal{C} \rightleftarrows L_{FE} \mathcal{D} : G$$

*is a symmetric monoidal Quillen equivalence.*

The major advantage of the Left Localization Principle over the Cellularization Principle is that the symmetric monoidal structure is preserved automatically. There are several variations of the Principle that we do not include in this introduction. Of particular note is the Compactly Generated Localization Principle, see Theorem 3.14. Although the assumptions of this last Principle are quite restrictive, there are interesting examples where it applies, as we show in our applications.

We now turn to the applications of the Left Localization Principle. The main motivation of the authors for developing the Left Localization Principle comes from rational equivariant stable homotopy theory.

**Algebraic models.** The programme of finding algebraic models for rational $G$-spectra was begun by Greenlees, who conjectured that for every compact Lie group $G$, there is an *abelian* category $\mathcal{A}(G)$, together with a Quillen equivalence between the category of rational $G$-spectra and the category of differential objects in $\mathcal{A}(G)$. The programme looks for abelian categories with finite homological dimension so that calculations can easily be performed, and equipped with an Adams spectral sequence to calculate homotopy classes of maps between $G$-spectra. This programme has so far been successful in the cases of $G$ finite [6], $G = SO(2)$ [44, 9] $G = O(2)$ [7], $G = SO(3)$ [33], $G$ a torus of any rank [28], the toral part of $G$-spectra [8], and free $G$-spectra for $G$ a compact Lie group [24, 26]. One can also ask for equivalences with extra structure such as being monoidal, so that the equivalence passes to ring and module spectra.

When attempting to find algebraic models for categories of interest, there are several techniques we can apply. One approach is to use Morita theory [43] which gives an equivalence with modules over the endomorphism ring of a generator. However, the endomorphism ring need not be commutative so that formality arguments are inaccessible, and the module category often has infinite homological dimension. Another alternative is to use the Cellularization Principle to reduce the problem to checking conditions on generating cells. In this paper, we show that the Left Localization Principle is another technique that we can use. Balchin-Greenlees [4] show that stable model categories can be split into pieces determined by left localizations in an adelic fashion, by proving that the stable model category is a homotopy pullback of an 'adelic cube'. We hope that the Left Localization Principle may be applied in these situations as well, to simplify the adelic cube.

**Completions.** In order to verify the conjecture of Greenlees in our case of interest, we discuss some homotopical aspects of completion. We briefly recall the relevant results about the different types of completions in algebra and we refer the reader to Section 5 for a more detailed exposition and references.

Let $I$ be a finitely generated ideal in a commutative ring $R$. The $I$-adic completion functor is a fundamental tool in algebra, but has poor homological properties as it is neither left nor right exact. Our approach is to work with its zeroth left derived functor which we denote by $L_0^I$. We say that an $R$-module $M$ is $L_0^I$-complete if the canonical map $M \to L_0^I M$ is an isomorphism. The full subcategory of $L_0^I$-complete modules is a symmetric monoidal abelian category which supports a projective model structure under a mild condition on the ideal considered. This condition is called weak pro-regularity and holds in many cases; for example, any ideal in a Noetherian ring is weakly pro-regular.

For homotopical purposes it is often convenient to consider the derived $I$-completion functor. This is defined in terms of the stable Koszul complex whose filtration provides a spectral sequence making the derived completion accessible. Under the weak pro-regularity hypothesis on the ideal $I$, the derived $I$-completion functor is equivalent to the total left derived functor of $I$-adic completion, and therefore calculates the local homology modules, see [21, 37].

We give a proof using the language of model categories that derived $I$-complete modules can be modelled via the abelian category of $L_0^I$-complete modules, see Theorem 6.10. It follows that a dg-module is derived $I$-complete if and only if its homology is $L_0^I$-complete. This generalises a result of Dwyer-Greenlees [17, 6.15] and clarifies an observation of Porta-Shaul-Yekutieli [37, 4.33] that derived $I$-complete modules need not have $I$-adically complete homology. We note the related work of Barthel-Heard-Valenzuela who have given an $\infty$-categorical approach to derived completion in the general setup of comodules over Hopf algebroids [12].

**Rational cofree $G$-spectra.** The equivariant stable homotopy category contains two classes of objects of particular note: the free and cofree $G$-spectra. An algebraic model for rational free $G$-spectra was constructed by Greenlees-Shipley [24, 26] in terms of torsion modules over the group cohomology ring. However, the abelian category of torsion modules is not monoidal as it has no tensor unit and therefore the Quillen equivalence in the free case cannot be refined to a symmetric monoidal Quillen equivalence.

By exploiting the equivalence between free and cofree $G$-spectra, we give a symmetric monoidal algebraic model for the category of rational cofree $G$-spectra. For convenience, we only state the result for the connected case in this introduction. See Theorem 9.6 for the general case.

**Theorem** (8.4)**.** *Let $G$ be a connected compact Lie group and $I$ be the augmentation ideal of $H^*BG$. Then there is a symmetric monoidal Quillen equivalence*

$$\mathrm{Sp}_G^{\mathrm{cofree}} \simeq_Q \mathrm{Mod}_{H^*BG}^{\wedge}$$

*between rational cofree $G$-spectra and $L_0^I$-complete dg-$H^*BG$-modules. In particular, there is a tensor-triangulated equivalence*

$$\text{cofree } G\text{-spectra} \simeq_\triangle \mathcal{D}(L_0^I\text{-complete } H^*BG\text{-modules}).$$

In this application, the Left Localization Principle manifests its advantages over the Cellularization Principle. Firstly, the proof of the equivalence is formal as it only requires a few elementary iterations of the Left Localization Principle and some formality arguments in algebra. In particular we avoid any "topological" formality argument using the Adams spectral sequence. Secondly, it gives a tensor-triangulated equivalence of the homotopy categories.

Free and cofree $G$-spectra are interesting for three particular reasons. Firstly, they represent cohomology theories on free $G$-spaces, the most prominent example of which is Borel cohomology. Secondly, the techniques employed in the construction of the algebraic models for free and cofree $G$-spectra are instructive for more general cases, such as that of torus-equivariant spectra [28]. Finally, the algebraic models for free and cofree $G$-spectra fit in the general picture of a local duality context in the sense of [11]. This means that the equivalence between free and cofree $G$-spectra in equivariant stable homotopy theory translates to the equivalence between torsion and complete modules in algebra.

**Contribution of this paper and related work.** Let us restrict to connected groups for simplicity, and continue to write $I$ for the augmentation ideal. A Quillen equivalence between rational cofree $G$-spectra and derived complete $H^*BG$-modules was already known by passing through free $G$-spectra in the following way:

$$\begin{array}{ccc} \text{free } G\text{-spectra} & \xleftrightarrow{\simeq_Q} & I\text{-power torsion-}H^*BG\text{-modules} \\ \simeq_Q \updownarrow & & \updownarrow \simeq_Q \\ \text{cofree } G\text{-spectra} & \longleftrightarrow & \text{derived } I\text{-complete-}H^*BG\text{-modules}. \end{array}$$

The horizontal Quillen equivalence is the algebraic model for free $G$-spectra of Greenlees-Shipley [24] and the right vertical follows from Dwyer-Greenlees' Morita theory [17] together with [25, §5]. However this is unsatisfactory for two main reasons. Firstly, it cannot be refined to a symmetric monoidal Quillen equivalence since the category of $I$-power torsion modules has no tensor unit. Secondly, it does not give an abelian model as desired in the conjecture of Greenlees. In light of this, our contribution is threefold: we prove the algebraic model for rational cofree $G$-spectra directly, we upgrade it to a symmetric monoidal Quillen equivalence, and we give an abelian model for derived complete modules. In addition, we collect several results about homotopical aspects of algebraic completions which we believe will be of independent interest.

Although our strategy is analogous to that employed by Greenlees-Shipley in the study of free $G$-spectra, the tools we use differ. In particular, the Left Localization Principle which we develop is a new and key ingredient in our proof.

**Outline of the paper.** The paper is divided into two main parts.

In the first part we give some necessary background on left Bousfield localizations and then state and prove the Left Localization Principle. We then investigate the implications in the case of homological localizations, which provide many key examples.

In the second part of the paper we focus on the applications of the Left Localization Principle. We apply the Left Localization Principle to understand completions of module categories and to construct a symmetric monoidal algebraic model for rational cofree $G$-spectra. We have decided to first construct the algebraic model for a connected compact Lie group and then show how to generalize our proofs to the non-connected case. In the final section, we construct a strongly convergent Adams spectral sequence to calculate homotopy classes of maps between cofree $G$-spectra.

**Conventions.** We shall follow the convention of writing the left adjoint above the right adjoint in an adjoint pair. We will use $q \colon QX \to X$ and $r \colon X \to RX$ to denote cofibrant and fibrant replacements of $X$ respectively.

**Acknowledgements.** We are extremely grateful to John Greenlees for many helpful discussions and suggestions. We would also like to thank Scott Balchin, Magdalena Kędziorek and Gabriel Valenzuela for their interest and comments.

**Part** 1. **The Left Localization Principle**

2. LEFT BOUSFIELD LOCALIZATION OF MODEL CATEGORIES

In this section we recall some necessary background on left Bousfield localizations following [29] and [10].

**Definition 2.1.** Let $\mathcal{C}$ be a model category and let $S$ be a collection of maps in $\mathcal{C}$.

- An object $W$ in $\mathcal{C}$ is *S-local* if it is fibrant in $\mathcal{C}$ and for every $s\colon A \to B$ in $S$, the natural map $\operatorname{Map}(B,W) \to \operatorname{Map}(A,W)$ is a weak equivalence of homotopy function complexes.
- A map $f\colon X \to Y$ in $\mathcal{C}$ is an *S-local equivalence* if for every $S$-local object $W$, the natural map $\operatorname{Map}(Y,W) \to \operatorname{Map}(X,W)$ is a weak equivalence of homotopy function complexes.

**Remark 2.2.** If the model category is stable, then the homotopy function complexes in the previous Definition can be replaced with the graded set of maps in the homotopy category, see [10, 4.5].

In many cases, we can endow the model category $\mathcal{C}$ with a new model structure, the *left Bousfield localization* of $\mathcal{C}$, in which the weak equivalences are the $S$-local equivalences, the cofibrations are unchanged, and the fibrant objects are the $S$-local objects. If it exists, we denote this model category by $L_S\mathcal{C}$.

**Hypothesis 2.3.** Throughout this paper we assume that all the required left Bousfield localizations exist.

**Remark 2.4.** The left Bousfield localization exists under mild conditions on the model category $\mathcal{C}$. For example, when $\mathcal{C}$ is left proper, cellular and $S$ is a set [29, 4.1.1], or when $\mathcal{C}$ is left proper, combinatorial and $S$ is a set [14, 4.7]. In particular, left Bousfield localizations (at sets of morphisms) exist for the stable model structure on spectra [35, 9.1], the stable model structure on equivariant spectra for any compact Lie group [34, III.4.2] and the projective model structure on dg-modules [13, 3.3].

Recall that a model category is *symmetric monoidal* if it is a closed symmetric monoidal category and it satisfies the *pushout-product axiom*: if $f\colon A \to B$ and $g\colon X \to Y$ are cofibrations, then the pushout-product map
$$f\square g\colon A \otimes Y \bigcup_{A \otimes X} B \otimes X \to B \otimes Y$$
is a cofibration, which is acyclic if either $f$ or $g$ is acyclic; and the *unit axiom*: the natural map $Q\mathbb{1} \otimes X \to \mathbb{1} \otimes X \cong X$ is a weak equivalence for all cofibrant $X$. We denote the internal hom functor by $F(-,-)$.

**Definition 2.5.** We say that a stable model category $\mathcal{C}$ is *homotopically compactly generated* by a set $\mathcal{K}$ of objects if its homotopy category $h\mathcal{C}$ is compactly generated by $\mathcal{K}$:

- for all $K \in \mathcal{K}$ and collections $\{M_i\}$ of objects of $\mathcal{C}$, the natural map $\bigoplus h\mathcal{C}(K, M_i) \to h\mathcal{C}(K, \bigoplus M_i)$ is an isomorphism;
- an object $X$ of $h\mathcal{C}$ is trivial if and only if $h\mathcal{C}(\Sigma^n K, X) = 0$ for all $K \in \mathcal{K}$ and $n \in \mathbb{Z}$.

Next we recall the definition of a monoidal Quillen adjunction from [42].

**Definition 2.6.** Let $F\colon \mathcal{C} \rightleftarrows \mathcal{D} : G$ be a Quillen adjunction between symmetric monoidal model categories.

(1) We say that $(F, G)$ is a *weak symmetric monoidal Quillen adjunction* if the right adjoint $G$ is lax monoidal (which gives the left adjoint $F$ an oplax monoidal structure) and the following two conditions hold:
   (a) for cofibrant $A$ and $B$ in $\mathcal{C}$, the oplax monoidal structure map $\phi\colon F(A \otimes B) \to F(A) \otimes F(B)$ is a weak equivalence in $\mathcal{D}$
   (b) for a cofibrant replacement $Q\mathbb{1}_\mathcal{C}$ of the unit in $\mathcal{C}$, the map $\phi_0\colon F(Q\mathbb{1}_\mathcal{C}) \to \mathbb{1}_\mathcal{D}$ is a weak equivalence in $\mathcal{D}$.
(2) If the oplax monoidal structure maps $\phi$ and $\phi_0$ are isomorphisms, then we say that $(F, G)$ is a *strong symmetric monoidal Quillen pair*.
(3) We say that the adjunction $(F, G)$ is *symmetric monoidal* if it is a weak symmetric monoidal Quillen adjunction.
(4) We say that the adjucntion $(F, G)$ is a *symmetric monoidal Quillen equivalence* if it is a symmetric monoidal adjunction and a Quillen equivalence.

**Remark 2.7.** A Quillen adjunction is symmetric monoidal if the left adjoint is strong monoidal and the unit object of $\mathcal{C}$ is cofibrant.

**Definition 2.8.** A set of morphisms $S$ of a stable model category $\mathcal{C}$ is said to be *stable* if the collection of $S$-local objects is closed under (de)suspensions. We say that a stable set of cofibrations $S$ of a stable, cellular, symmetric monoidal model category $\mathcal{C}$ is *monoidal* if $S\square I = \{s\square i \mid s \in S,\ i \in I\}$ is contained in the class of $S$-equivalences, where $I$ is the set of generating cofibrations for $\mathcal{C}$.

We will need the following result.

**Proposition 2.9** ([10, 5.1]). *Let $\mathcal{C}$ be a proper, cellular, stable, symmetric monoidal model category and let $S$ be a stable set of cofibrations between cofibrant objects. Then the localization $L_S\mathcal{C}$ is a symmetric monoidal model category if and only if $S$ is monoidal.*

**Remark 2.10.** Any map in a model category can be replaced up to weak equivalence by a cofibration between cofibrant objects: first cofibrantly replacing the source and then factoring the composite into a cofibration followed by an acyclic fibration. Since left Bousfield localization depends only on the homotopy type of the class of maps, we can assume without loss of generality that $S$ consists of cofibrations between cofibrant objects.

## 3. The Left Localization Principle

We are now ready to work towards the Left Localization Principle. Before we can prove an induced Quillen equivalence, we must check that the Quillen adjunction descends to the localizations. Recall that $Q$ and $R$ denote cofibrant and fibrant replacement in the original model structures on $\mathcal{C}$ and $\mathcal{D}$ respectively.

**Proposition 3.1.** *Let $F : \mathcal{C} \rightleftarrows \mathcal{D} : G$ be a Quillen adjunction between stable model categories satisfying Hypothesis 2.3.*

(1) *Let $S$ be a stable set of morphisms in $\mathcal{C}$. Then the adjunction*

$$F : L_S\mathcal{C} \rightleftarrows L_{FQS}\mathcal{D} : G$$

*is a Quillen adjunction. Furthermore, it is a symmetric monoidal Quillen adjunction if $F : \mathcal{C} \rightleftarrows \mathcal{D} : G$ is a symmetric monoidal Quillen adjunction and $S$ and $FQS$ are monoidal.*

(2) *Let $T$ be a stable set of morphisms in $\mathcal{D}$. Suppose that $F$ sends $GRT$-equivalences between cofibrant objects to $T$-equivalences. Then the adjunction*

$$F : L_{GRT}\mathcal{C} \rightleftarrows L_T\mathcal{D} : G$$

*is a Quillen adjunction. Furthermore, it is a symmetric monoidal Quillen adjunction if $F : \mathcal{C} \rightleftarrows \mathcal{D} : G$ is a symmetric monoidal Quillen adjunction and $T$ and $GRT$ are monoidal.*

*Proof.* Let us prove (1). By Dugger [16, A.2] it is sufficient to check that $G$ preserves fibrations between fibrant objects and all acyclic fibrations. The acyclic fibrations in a left Bousfield localization are the same as in the original model structure so $G$ preserves them. Since fibrations between local objects in a left Bousfield localization are the same as original fibrations between local objects, it is sufficient to check that $G$ sends $FQS$-local objects to $S$-local objects. Suppose that $W$ is $FQS$-local. Let $s : A \to B$ be in $S$. We have

$$\begin{array}{ccc} \mathrm{Map}(B, GW) & \longrightarrow & \mathrm{Map}(A, GW) \\ \sim\downarrow & & \downarrow\sim \\ \mathrm{Map}(FQB, W) & \underset{\sim}{\longrightarrow} & \mathrm{Map}(FQA, W) \end{array}$$

so that $GW$ is $S$-local.

The claim about the monoidality follows from the fact that the cofibrations in a left Bousfield localization are the same as in the original category, and the local equivalences contain the original weak equivalences.

Let us now prove (2). By Hirschhorn [29, 3.3.18], to prove that $L_{GRT}\mathcal{C} \rightleftarrows L_T\mathcal{D}$ is a Quillen adjunction, it is sufficient to check that $F$ sends $GRT$-equivalences between cofibrant objects to $T$-equivalences which was a hypothesis. The claim about monoidal structures is clear. $\square$

**Remark 3.2.** We note that the hypothesis that $F$ sends $GRT$-equivalences between cofibrant objects to $T$-equivalences in part (2) of the previous Proposition may seem hard to verify in practice. However, we show in Lemma 3.12 that in the case of homological localization, this hypothesis can be replaced by a condition which is much easier to verify.

**Remark 3.3.** If $S$ is monoidal, it often happens that $FQS$ is also monoidal. Write $I_\mathcal{C}$ and $I_\mathcal{D}$ for the sets of generating cofibrations in $\mathcal{C}$ and $\mathcal{D}$ respectively. For instance, one can easily check that $FQS$ is monoidal when $(F,G)$ is a strong symmetric monoidal Quillen pair and $I_\mathcal{D} \subseteq F(I_\mathcal{C})$, or, when $(F,G)$ is a weak symmetric monoidal Quillen pair, the domains of $I_\mathcal{C}$ are cofibrant and $I_\mathcal{D} \subseteq F(I_\mathcal{C})$. Note that the condition that $I_\mathcal{D} \subseteq F(I_\mathcal{C})$ is satisfied in the case when the model structure on $\mathcal{D}$ is right induced from $\mathcal{C}$.

We can now state and prove the Left Localization Principle. We note that as the cofibrations are the same in the left Bousfield localization as in the original model structure, we continue to write $Q$ for the cofibrant replacement in the localization. However, since being fibrant in the localization is a stronger condition than being fibrant in the original model structure, we write $\overline{R}$ for the fibrant replacement in the localization.

**Theorem 3.4** (Left Localization Principle)**.** *Let $\mathcal{C}$ and $\mathcal{D}$ be stable model categories satisfying Hypothesis 2.3 and let $F : \mathcal{C} \rightleftarrows \mathcal{D} : G$ be a Quillen adjunction.*

(1) *Suppose that $\mathcal{C}$ is homotopically compactly generated by a set $\mathcal{K}$ and that $\mathcal{D}$ is homotopically compactly generated by $FQ\mathcal{K}$. Let $S$ be a stable set of morphisms in $\mathcal{C}$. Suppose that the following conditions hold:*
  (i) *The derived unit map $\eta_K \colon QK \to GRFQK$ is an $S$-equivalence for all $K \in \mathcal{K}$;*
  (ii) *$G$ sends $FQS$-equivalences between fibrant objects in $\mathcal{D}$ to $S$-equivalences.*
  *Then the induced Quillen adjunction*

$$F : L_S\mathcal{C} \rightleftarrows L_{FQS}\mathcal{D} : G$$

  *is a Quillen equivalence. Moreover, if $F : \mathcal{C} \rightleftarrows \mathcal{D} : G$ is a symmetric monoidal Quillen adjunction and $S$ and $FQS$ are monoidal, then $F : L_S\mathcal{C} \rightleftarrows L_{FQS}\mathcal{D} : G$ is a symmetric monoidal Quillen equivalence.*

(2) *Suppose that $\mathcal{D}$ is homotopically compactly generated by a set $\mathcal{L}$ and that $\mathcal{C}$ is homotopically compactly generated by $GR\mathcal{L}$. Let $T$ be a stable set of morphisms in $\mathcal{D}$. Suppose that the following conditions hold:*
  (i) *The derived counit map $\epsilon_L \colon FQGRL \to RL$ is a $T$-equivalence for all $L \in \mathcal{L}$[1];*
  (ii) *$G$ sends $T$-equivalences between fibrant objects in $\mathcal{D}$ to $GRT$-equivalences.*
  (iii) *$F$ sends $GRT$-equivalences between cofibrant objects to $T$-equivalences;*
  *Then the induced Quillen adjunction*

$$F : L_{GRT}\mathcal{C} \rightleftarrows L_T\mathcal{D} : G$$

  *is a Quillen equivalence. Moreover, if $F : \mathcal{C} \rightleftarrows \mathcal{D} : G$ is a symmetric monoidal Quillen adjunction and $T$ and $GRT$ are monoidal, then $F : L_{GRT}\mathcal{C} \rightleftarrows L_T\mathcal{D} : G$ is a symmetric monoidal Quillen equivalence.*

*Proof.* Let us prove (1). Firstly we show that the derived functor $GR$ preserves sums, so that the subcategories

$$\mathcal{A} = \{X \in \mathrm{h}\mathcal{C} \mid \eta_X \colon QX \xrightarrow{\sim_S} GRFQX\} \quad \text{and} \quad \mathcal{A}' = \{Y \in \mathrm{h}\mathcal{D} \mid \epsilon_Y \colon FQGRY \xrightarrow{\sim_{FQS}} RY\}$$

are localizing. Let $(X_i)_{i \in I}$ be a collection of objects in $\mathrm{h}\mathcal{D}$. Using compactness we see that for all $K \in \mathcal{K}$

$$\mathrm{h}\mathcal{C}(K, GR(\bigoplus_{i \in I} X_i)) \cong \mathrm{h}\mathcal{D}(FQK, \bigoplus_{i \in I} X_i) \cong \bigoplus_{i \in I} \mathrm{h}\mathcal{D}(FQK, X_i) \cong \bigoplus_{i \in I} \mathrm{h}\mathcal{C}(K, GR(X_i)) \cong \mathrm{h}\mathcal{C}(K, \bigoplus_{i \in I} GR(X_i)).$$

Since $\mathcal{K}$ generates $\mathrm{h}\mathcal{C}$ we conclude that $GR$ preserves arbitrary sums.

---
[1]This condition should instead read: "The derived counit map $\epsilon_L \colon FQGRL \to RL$ is a weak equivalence in $\mathcal{D}$ for all $L \in \mathcal{L}$"; see the corrigendum at the end of this paper

By assumption (i), we know that $\mathcal{K} \subset \mathcal{A}$ thus $\mathcal{A} = h\mathcal{C}$ as $\mathcal{K}$ generates $h\mathcal{C}$. Note that $FQ\eta_K$ is an $FQS$-equivalence. Using the triangular identity of the derived adjunction

$$FQK \xrightarrow{FQ\eta_K} FQGRFQK$$
$$\searrow_r \quad \downarrow_{\epsilon_{FQK}}$$
$$RFQK$$

and 2-out-of-3, we obtain that $FQK \in \mathcal{A}'$ and hence $\mathcal{A}' = h\mathcal{D}$ as $FQ\mathcal{K}$ generates $h\mathcal{D}$.

We must prove that $\overline{\eta}_X \colon QX \to G\overline{R}FQX$ is an $S$-equivalence for all $X \in h\mathcal{C}$ and that $\overline{\epsilon}_Y \colon FQG\overline{R}Y \to \overline{R}Y$ is an $FQS$-equivalence for all $Y \in h\mathcal{D}$. Note that the canonical map $GRFQX \to G\overline{R}FQX$ is an $S$-equivalence by condition (ii). Therefore the derived unit

$$\overline{\eta}_X \colon QX \xrightarrow{\sim_S} GRFQX \xrightarrow{\sim_S} G\overline{R}FQX$$

is an $S$-equivalence. For the derived counit, note that the canonical map $GRY \to G\overline{R}Y$ is an $S$-equivalence and therefore $FQGRY \to FQG\overline{R}Y$ is an $FQS$-equivalence by Ken Brown's Lemma. By considering the diagram

$$\begin{array}{ccc} FQGRY & \xrightarrow{\epsilon_Y} & RY \\ \sim_{FQS} \downarrow & & \downarrow \sim_{FQS} \\ FQG\overline{R}Y & \xrightarrow{\overline{\epsilon}_Y} & \overline{R}Y \end{array}$$

we see that $\overline{\epsilon}_Y$ is an $FQS$-equivalence if and only if $\epsilon_Y$ is so. Since $\mathcal{A}' = h\mathcal{D}$ the claim follows.

The proof of part (2) is essentially the same. Note that condition (iii) comes into play only to ensure that the Quillen adjunction descends to the localizations, see Proposition 3.1. □

**Remark 3.5.** Notice that the conditions in (1) imply that the derived functor $FQ$ preserves all compact objects. Moreover, in the proof we showed that $GR$ preserves sums so it also follows that under the conditions in (2) the derived functor $GR$ preserves all compact objects.

**Remark 3.6.** In [30, 2.3] Hovey gives criteria for when left Bousfield localization preserves Quillen equivalences. His result does not assume stability but does not treat the case where the original adjunction is not a Quillen equivalence.

In the Left Localization Principle we assumed that $\mathcal{C}$ and $\mathcal{D}$ are homotopically compactly generated whereas in the following we assume that the localizations are homotopically compactly generated. This is a stronger condition but holds in certain cases when the localization is homological, see Remark 3.16.

**Theorem 3.7** (Compactly Generated Localization Principle)**.** *Let $\mathcal{C}$ and $\mathcal{D}$ be stable model categories satisfying Hypothesis 2.3 and let $F \colon \mathcal{C} \rightleftarrows \mathcal{D} \colon G$ be a Quillen adjunction.*

(1) *Let $S$ be a stable set of morphisms in $\mathcal{C}$. Suppose that $L_S\mathcal{C}$ is homotopically compactly generated by a set $\mathcal{K}$ and that $L_{FQS}\mathcal{D}$ is homotopically compactly generated by $FQ\mathcal{K}$. If the derived unit map $\overline{\eta}_K \colon QK \to G\overline{R}FQK$ is an $S$-equivalence for all $K \in \mathcal{K}$ then the induced Quillen adjunction*

$$F \colon L_S\mathcal{C} \rightleftarrows L_{FQS}\mathcal{D} \colon G$$

*is a Quillen equivalence. Moreover, if $F \colon \mathcal{C} \rightleftarrows \mathcal{D} \colon G$ is a symmetric monoidal Quillen adjunction and $S$ and $FQS$ are monoidal, then $F \colon L_S\mathcal{C} \rightleftarrows L_{FQS}\mathcal{D} \colon G$ is a symmetric monoidal Quillen equivalence.*

(2) *Let $T$ be a stable set of morphisms in $\mathcal{D}$. Suppose that $L_T\mathcal{D}$ is homotopically compactly generated by a set $\mathcal{L}$ and that $L_{GRT}\mathcal{C}$ is homotopically compactly generated by $GR\mathcal{L}$. Suppose that the derived counit $\overline{\epsilon}_L \colon FQG\overline{R}L \to \overline{R}L$ is a $T$-equivalence for all $L \in \mathcal{L}^2$ and that $F$ sends $GRT$-equivalences between cofibrant objects to $T$-equivalences. Then the induced Quillen adjunction*

$$F \colon L_{GRT}\mathcal{C} \rightleftarrows L_T\mathcal{D} \colon G$$

---

[2]This condition should instead read: "Suppose that the derived counit map $\epsilon_L \colon FQGRL \to RL$ is a weak equivalence in $\mathcal{D}$ for all $L \in \mathcal{L}$"; see the corrigendum at the end of this paper

*is a Quillen equivalence. Moreover, if $F : \mathcal{C} \rightleftarrows \mathcal{D} : G$ is a symmetric monoidal Quillen adjunction and $T$ and $GRT$ are monoidal, then $F : L_{GRT}\mathcal{C} \rightleftarrows L_T\mathcal{D} : G$ is a symmetric monoidal Quillen equivalence.*

*Proof.* Apply the Cellularization Principle [25, 2.7] to the Quillen adjunctions $F : L_S\mathcal{C} \rightleftarrows L_{FQS}\mathcal{D} : G$ and $F : L_{GRT}\mathcal{C} \rightleftarrows L_T\mathcal{D} : G$ obtained from Proposition 3.1. □

3.1. **Homological localization.** We now rephrase the Left Localization Principle for homological Bousfield localizations. This setting provides a large family of examples in which our result simplifies.

**Definition 3.8.** Let $\mathcal{C}$ be a stable and symmetric monoidal model category and let $E$ be a cofibrant object of $\mathcal{C}$. We say that $f : X \to Y$ is an *E-equivalence* if $E \otimes f : E \otimes X \to E \otimes Y$ is a weak equivalence.

When it exists, localizing at the $E$-equivalences produces a model structure on $\mathcal{C}$ in which the weak equivalences are the $E$-equivalences, the cofibrations are unchanged and the fibrant objects are the $E$-local objects. We call this new model category the *homological localization* of $\mathcal{C}$ at $E$ and write $L_E\mathcal{C}$.

**Hypothesis 3.9.** From now on we assume that the required homological localizations exist.

**Remark 3.10.** The homological localization exists if $\mathcal{C}$ is a stable, symmetric monoidal, proper and compactly generated model category in the sense of [46, 1.2.3.4]; see [18, §VIII.1] for the special case of spectra, and [4, §6.A] for the general case.

**Proposition 3.11.** *Let $\mathcal{C}$ be a symmetric monoidal model category satisfying Hypothesis 3.9, and let $E$ be a cofibrant object of $\mathcal{C}$. Then the homological localization $L_E\mathcal{C}$ is a symmetric monoidal model category.*

*Proof.* Take two cofibrations $i$ and $j$. Since the cofibrations in $L_E\mathcal{C}$ are the same as in $\mathcal{C}$, the pushout-product map $i\square j$ is a cofibration since $\mathcal{C}$ satisfies the pushout-product axiom. Now suppose that $i$ is an $E$-equivalence also. We have that $E\otimes(i\square j) = (E\otimes i)\square(E\otimes j)$ since $E\otimes-$ is a left adjoint. The functor $E\otimes-$ is left Quillen since $E$ is cofibrant, so $E \otimes i$ is an acyclic cofibration and $E \otimes j$ is a cofibration. Therefore, $E \otimes (i\square j)$ is a weak equivalence by the pushout-product axiom for $\mathcal{C}$. In other words, $i\square j$ is an $E$-equivalence as required. The unit axiom follows immediately from the unit axiom for $\mathcal{C}$, since the cofibrations are the same in the left Bousfield localization as in the original model category. □

**Lemma 3.12.** *Let $F : \mathcal{C} \rightleftarrows \mathcal{D} : G$ be a symmetric monoidal Quillen adjunction between stable symmetric monoidal model categories and let $E'$ be a bifibrant object in $\mathcal{D}$. If $\epsilon_E : FQGE' \to E'$ is a weak equivalence in $\mathcal{D}$, then $F$ sends $QGE'$-equivalences between cofibrant objects to $E'$-equivalences.*

*Proof.* Let $X \to Y$ be a $QGE'$-equivalence between cofibrant objects. By Ken Brown's lemma, $F(QGE' \otimes X) \to F(QGE' \otimes Y)$ is a weak equivalence. We have the commutative diagram

$$\begin{array}{ccccc} F(QGE' \otimes X) & \xrightarrow{\sim} & FQGE' \otimes FX & \xrightarrow{\sim} & E' \otimes FX \\ \sim\downarrow & & \downarrow & & \downarrow \\ F(QGE' \otimes Y) & \xrightarrow{\sim} & FQGE' \otimes FY & \xrightarrow{\sim} & E' \otimes FY \end{array}$$

in which the first horizontal maps are equivalences by definition of a symmetric monoidal Quillen pair, and the second horizontal maps are equivalences since $\epsilon_E : FQGE' \to E'$ is a weak equivalence and tensoring with a cofibrant object preserves weak equivalences between cofibrants by Ken Brown's lemma. Hence by two-out-of-three, $E' \otimes FX \to E' \otimes FY$ is a weak equivalence as required. □

Recall that the homological localization at an object $E$ is a special case of left Bousfield localization which inverts the $E$-equivalences. Therefore we can combine this Lemma with the Left Localization Principle to obtain our version for homological localizations.

**Theorem 3.13** (Left Localization Principle). *Let $\mathcal{C}$ and $\mathcal{D}$ be stable, symmetric monoidal model categories satisfying Hypothesis 3.9 and let $F : \mathcal{C} \rightleftarrows \mathcal{D} : G$ be a symmetric monoidal Quillen adjunction.*

(1) *Suppose that $\mathcal{C}$ is homotopically compactly generated by a set $\mathcal{K}$ and that $\mathcal{D}$ is homotopically compactly generated by $FQ\mathcal{K}$. Let $E \in \mathcal{C}$ be cofibrant. Suppose that the following conditions hold:*
   (i) *The derived unit map $QK \to GRFQK$ is an $E$-equivalence for all $K \in \mathcal{K}$;*
   (ii) *$G$ sends $FE$-equivalences between fibrant objects in $\mathcal{D}$ to $E$-equivalences.*
   *Then the induced Quillen adjunction*
   $$F: L_E\mathcal{C} \rightleftarrows L_{FE}\mathcal{D} : G$$
   *is a symmetric monoidal Quillen equivalence.*

(2) *Suppose that $\mathcal{D}$ is homotopically compactly generated by a set $\mathcal{L}$ and that $\mathcal{C}$ is homotopically compactly generated by $GR\mathcal{L}$. Let $E' \in \mathcal{D}$ be bifibrant. Suppose that the following conditions hold:*
   (i) *The derived counit map $FQGRL \to RL$ is an $E'$-equivalence for all $L \in \mathcal{L}$;*
   (ii) *$G$ sends $E'$-equivalences between fibrant objects in $\mathcal{D}$ to $QGE'$-equivalences.*
   (iii) *The map $FQGE' \to E'$ is a weak equivalence in $\mathcal{D}$;*
   *Then the induced Quillen adjunction*
   $$F: L_{QGE'}\mathcal{C} \rightleftarrows L_{E'}\mathcal{D} : G$$
   *is a symmetric monoidal Quillen equivalence.*

We now give a mixing of the Left Localization Principle and the Cellularization Principle. Note that we again write $\overline{R}$ for a fibrant replacement in the Bousfield localization.

**Theorem 3.14** (Compactly Generated Localization Principle). *Let $\mathcal{C}$ and $\mathcal{D}$ be stable, symmetric monoidal model categories satisfying Hypothesis 3.9 and let $F: \mathcal{C} \rightleftarrows \mathcal{D} : G$ be a symmetric monoidal Quillen adjunction.*

(1) *Let $E$ be a cofibrant object of $\mathcal{C}$. Suppose that $L_E\mathcal{C}$ is homotopically compactly generated by a set $\mathcal{K}$ and that $L_{FE}\mathcal{D}$ is homotopically compactly generated by $FQ\mathcal{K}$. If the derived unit map $Q\overline{\eta}_K: K \to G\overline{R}FQK$ is an $E$-equivalence for all $K \in \mathcal{K}$ then the induced Quillen adjunction*
   $$F: L_E\mathcal{C} \rightleftarrows L_{FE}\mathcal{D} : G$$
   *is a symmetric monoidal Quillen equivalence.*

(2) *Let $E'$ be a bifibrant object of $\mathcal{D}$. Suppose that $L_{E'}\mathcal{D}$ is homotopically compactly generated by a set $\mathcal{L}$ and that $L_{QGE'}\mathcal{C}$ is homotopically compactly generated by $GR\mathcal{L}$. Suppose that the derived counit $\overline{\epsilon}_L: FQG\overline{R}L \to \overline{R}L$ is an $E'$-equivalence for all $L \in \mathcal{L}$ and that $FQG\overline{R}E' \to \overline{R}E'$ is a weak equivalence in $\mathcal{D}$. Then the induced Quillen adjunction*
   $$F: L_{QGE'}\mathcal{C} \rightleftarrows L_{E'}\mathcal{D} : G$$
   *is a symmetric monoidal Quillen equivalence.*

**Remark 3.15.** Barnes-Roitzheim have compared left and right Bousfield localizations of stable model categories at dualizable objects [10, 9.6]. More precisely, they proved that the identity functors
$$L_A\mathcal{C} \leftrightarrows \mathrm{Cell}_{DA}\mathcal{C}$$
give a Quillen equivalence, where $D = F(-, \mathbb{1})$ is the dual functor and $A$ is dualizable. Accordingly, in some cases the Left Localization Principle can be replaced by the Cellularization Principle and vice versa. However, there are some subtleties that need to be considered. Firstly, the two principles are "exchangeable" only if the functors interact well with taking duals and we localize at dualizable objects. This a big disadvantage for instance in global stable homotopy theory where almost no compact objects are dualizable. This was one of the main motivations of the authors to develop the Left Localization Principle. Secondly, the two principles have quite different behaviour when we take into account the symmetric monoidal structure. While the Left Localization Principle for homological localization automatically yields a monoidal Quillen equivalence, the Cellularization Principle requires strong conditions, in particular when passing cells along the right adjoint, see [9, 5.1.7].

**Remark 3.16.** If we want to apply the Compactly Generated Localization Principle we need to know that the category of local objects is compactly generated. This holds for instance, when we localize at dualizable objects. More precisely, let $\mathcal{C}$ be a stable, symmetric monoidal model category, and let $A$ be a dualizable object of $\mathcal{C}$. It is not difficult to see that if $\mathcal{C}$ is homotopically compactly generated by a set $\mathcal{K}$ then the homological localization $L_A\mathcal{C}$ is homotopically compactly generated by $DA \otimes \mathcal{K}$. Firstly, $DA \otimes K$ is $A$-local for all $K \in \mathcal{K}$ since if $A \otimes Z \simeq 0$, then $h\mathcal{C}(Z, DA \otimes K) = h\mathcal{C}(Z, F(A, K)) = h\mathcal{C}(Z \otimes A, K) = 0$. Compactness follows from the fact that $A \otimes - : hL_A\mathcal{C} \to h\mathcal{C}$ preserves colimits, and the generation is an immediate consequence of the duality adjunction. For more details, see for instance [36, 2.27].

## 4. Completion of module categories

In this section we apply the Left Localization Principle to obtain symmetric monoidal Quillen equivalences relating a ring to its completion. We provide a general statement and then give several concrete examples of interest.

**Notation 4.1.** Given a commutative monoid $R$ in a symmetric monoidal model category $\mathcal{C}$, we denote by $\mathrm{Mod}_R(\mathcal{C})$ the category of $R$-modules equipped with the projective model structure (if it exists) in which the weak equivalences and fibrations are created by the forgetful functor $\mathrm{Mod}_R(\mathcal{C}) \to \mathcal{C}$. If the underlying category is clear, we will often write $\mathrm{Mod}_R$.

**Hypothesis 4.2.** Throughout this paper we assume that the projective model structure on $\mathrm{Mod}_R(\mathcal{C})$ exists and that it is left proper, so that left Bousfield localizations exist.

**Remark 4.3.** Note that the projective model structure exists if $\mathcal{C}$ satisfies the monoid axiom [41, 4.1], and it is left proper in many cases: for instance in categories of (equivariant) spectra [35, 12.1(i)] and [34, III.7.6], and in dg-modules [13, 3.3].

**Proposition 4.4.** *Let $\mathcal{C}$ be a stable, symmetric monoidal model category, homotopically compactly generated by a set $\mathcal{K}$. Let $R$ be a commutative monoid in $\mathcal{C}$, let $E$ be a cofibrant $R$-module and suppose that the natural map $\theta \colon R \to L_E R$ is a map of commutative monoids. The natural map $\theta$ induces a symmetric monoidal extension-restriction of scalars Quillen adjunction*

$$L_E R \otimes_R - : \mathrm{Mod}_R(\mathcal{C}) \rightleftarrows \mathrm{Mod}_{L_E R}(\mathcal{C}) : \theta^*$$

*between the categories of modules. Then the Left Localization Principle applies and gives a symmetric monoidal Quillen equivalence*

$$L_E \mathrm{Mod}_R(\mathcal{C}) \simeq_Q L_E \mathrm{Mod}_{L_E R}(\mathcal{C}).$$

**Remark 4.5.** Note that there is an abuse of notation in the Proposition above since in general there is no natural $L_E R$-module structure on $E$ at the model category level. More precisely, on the right hand side of the Quillen equivalence above we should have localized at $L_E R \otimes_R E$ instead of $E$. However, this abuse of notation does no harm since there is a natural weak equivalence $E \xrightarrow{\sim} L_E R \otimes_R E$ in $\mathcal{C}$ and the class of $L_E R \otimes_R E$-equivalences is detected in the homotopy category of $\mathcal{C}$.

*Proof.* The set $R \otimes \mathcal{K}$ provides a set of compact generators for $h\mathrm{Mod}_R(\mathcal{C})$. The left adjoint is strong monoidal and maps compact generators to compact generators since $L_E R \otimes_R (R \otimes \mathcal{K}) \cong L_E R \otimes \mathcal{K}$.

Without loss of generality we may now assume that $\mathcal{K}$ consists of cofibrant objects. As $R \to L_E R$ is an $E$-equivalence, we obtain a weak equivalence $E \xrightarrow{\sim} L_E R \otimes_R E$ by tensoring with $E$. Therefore, the derived unit $R \otimes K \to L_E R \otimes_R (R \otimes K) = L_E R \otimes K$ is an $E$-equivalence for all $K \in \mathcal{K}$.

Finally we must show that the right adjoint $\theta^*$ preserves $E$-equivalences between fibrant objects. Note that there is a natural map $E \otimes_R \theta^* M \to \theta^*(E \otimes_{L_E R} M)$ of $R$-modules, which is a weak equivalence as $E \simeq L_E R \otimes_R E$. Now suppose that $M \to N$ is an $E$-equivalence between fibrant $L_E R$-modules. By considering the diagram

$$\begin{array}{ccc} E \otimes_R \theta^* M & \longrightarrow & E \otimes_R \theta^* N \\ \sim \downarrow & & \downarrow \sim \\ \theta^*(E \otimes_{L_E R} M) & \longrightarrow & \theta^*(E \otimes_{L_E R} N) \end{array}$$

we see that $\theta^* M \to \theta^* N$ is an $E$-equivalence of $R$-modules. $\square$

**Example 4.6.** Let $\mathbb{Z}_p$ denote the $p$-adic integers and consider the ring map $\theta \colon \mathbb{Z} \to \mathbb{Z}_p$ which induces a symmetric monoidal Quillen adjunction between the categories of chain complexes

$$\mathbb{Z}_p \otimes_{\mathbb{Z}} - \colon \mathrm{Mod}_{\mathbb{Z}} \rightleftarrows \mathrm{Mod}_{\mathbb{Z}_p} \colon \theta^*$$

via extension and restriction of scalars. We can apply Proposition 4.4 to obtain a symmetric monoidal Quillen equivalence

$$L_{\mathbb{Z}/p}\mathrm{Mod}_{\mathbb{Z}} \simeq_Q L_{\mathbb{Z}/p}\mathrm{Mod}_{\mathbb{Z}_p}.$$

By [22, 4.2], we can identify the homotopy categories of the two localizations with the subcategories of the derived categories consisting of derived $p$-complete modules which we denote $\Lambda_{\mathbb{Z}/p}\mathrm{Mod}_{\mathbb{Z}}$ and $\Lambda_{\mathbb{Z}/p}\mathrm{Mod}_{\mathbb{Z}_p}$ respectively. Putting everything together we get a tensor-triangulated equivalence

$$\Lambda_{\mathbb{Z}/p}\mathrm{Mod}_{\mathbb{Z}} \simeq_{\triangle} \Lambda_{\mathbb{Z}/p}\mathrm{Mod}_{\mathbb{Z}_p}.$$

**Example 4.7.** Let $G$ be a compact Lie group and $\mathcal{F}$ a family of subgroups of $G$. Note that the $G$-spectrum $DE\mathcal{F}_+$ is a commutative ring $G$-spectrum via the diagonal map $\Delta \colon E\mathcal{F}_+ \to E\mathcal{F}_+ \wedge E\mathcal{F}_+$. It is easy to check that $DE\mathcal{F}_+$ is $E\mathcal{F}_+$-local and that the unit map $\eta \colon S^0 \to DE\mathcal{F}_+$ is a $E\mathcal{F}_+$-equivalence. We can then apply Proposition 4.4 to obtain a symmetric monoidal Quillen equivalence

$$L_{E\mathcal{F}_+}\mathrm{Sp}_G \simeq_Q L_{E\mathcal{F}_+}\mathrm{Mod}_{DE\mathcal{F}_+}(\mathrm{Sp}_G).$$

Note that the proof of Proposition 4.4 works more generally for localizations at a set of maps $S$, provided that the natural map $\theta \colon R \to L_S R$ is a map of commutative monoids.

**Example 4.8.** Let $\mathcal{G}$ be the global family of compact Lie groups. Denote by $\mathrm{Sp}_{\mathcal{G}}$ the category of orthogonal spectra with the $\mathcal{G}$-global model structure which is proper [40, 4.3.17]. By [40, 4.5.21, 4.5.22(ii)], there exists a morphism of ultracommutative ring spectra $i_{\mathbb{S}} \colon \mathbb{S} \to \mathrm{b}\mathbb{S}$ between the global sphere spectrum and the global Borel construction which exhibits $\mathrm{b}\mathbb{S}$ as a localization of the global sphere spectrum at the class of non-equivariant equivalences. Note that the projective model structure on $\mathrm{Mod}_{\mathrm{b}\mathbb{S}}(\mathrm{Sp}_{\mathcal{G}})$ exists by [40, 4.3.29] and it is proper by a similar argument as in [35, 12.1(i)] so that we can perform left Bousfield localizations. Following the proof of Proposition 4.4 and localizing at the class 1 of non-equivariant equivalences (see Remark 4.9 for justification of its existence), we obtain a symmetric monoidal Quillen equivalence

$$L_1 \mathrm{Sp}_{\mathcal{G}} \simeq_Q L_1 \mathrm{Mod}_{\mathrm{b}\mathbb{S}}(\mathrm{Sp}_{\mathcal{G}}).$$

We note that this is a symmetric monoidal Quillen equivalence using Remark 3.3, since the model structure on $\mathrm{Mod}_{\mathrm{b}\mathbb{S}}(\mathrm{Sp}_{\mathcal{G}})$ is right induced from the $\mathcal{G}$-global model structure on $\mathrm{Sp}_{\mathcal{G}}$. Finally using the language of [40] we can identify the homotopy category of $L_1\mathrm{Sp}_{\mathcal{G}}$ with the full subcategory of the global stable homotopy category consisting of those global spectra which are right induced from the trivial family.

**Remark 4.9.** It is not immediate that the left Bousfield localization of $\mathrm{Sp}_{\mathcal{G}}$ at the *class* of non-equivariant equivalence actually exists. This localization cannot be constructed as a homological localization since in global stable homotopy theory an analogue of the free $G$-space $EG$ does not exist. Instead we apply Bousfield-Friedlander localization [15, 9.3] to the natural transformation $i_X \colon X \to \mathrm{b}X$ which is a non-equivariant equivalence. By construction, the global Borel functor b has the property that for all $G \in \mathcal{G}$, the underlying $G$-spectrum of $\mathrm{b}X$ is cofree, see [40, 4.5.16, 4.5.22]. In particular this shows that $f \colon X \to Y$ is a non-equivariant equivalence if and only if $\mathrm{b}f \colon \mathrm{b}X \to \mathrm{b}Y$ is a global equivalence. The conditions (A1) and (A2) from [15, 9.2] easily follow from this observation. The final condition (A3) follows from the right properness of $\mathrm{Sp}_{\mathcal{F}}$ for the trivial family $\mathcal{F} = \{1\}$, together with the fact that any $\mathcal{G}$-global fibration is a $\mathcal{F}$-global fibration. The argument for $\mathrm{Mod}_{\mathrm{b}\mathbb{S}}(\mathrm{Sp}_{\mathcal{G}})$ is similar as the weak equivalences and fibrations are created by the forgetful functor $\mathrm{Mod}_{\mathrm{b}\mathbb{S}}(\mathrm{Sp}_{\mathcal{G}}) \to \mathrm{Sp}_{\mathcal{G}}$.

## Part 2. Rational cofree $G$-spectra

We give a symmetric monoidal algebraic model for the category of rational cofree $G$-spectra for $G$ a compact Lie group, in the sense of [19]. We will initially prove the result for $G$ connected and then show how to extend our proofs to any compact Lie group. In the final section we construct a strongly convergent Adams spectral sequence calculating homotopy classes of maps between cofree $G$-spectra.

## 5. Completions in algebra

We now recall some results about complete modules following [21].

Let $R$ be a graded commutative ring and let $I$ be a finitely generated homogeneous ideal. The $I$-adic completion of a module $M$ is defined by
$$M_I^\wedge = \lim_n M/I^n M.$$
We say that a module $M$ is *$I$-adically complete* if the natural map $M \to M_I^\wedge$ is an isomorphism. A dg-module is said to be *$I$-adically complete* if its underlying graded module is.

Since the $I$-adic completion functor is neither left nor right exact in general, our approach is to consider the zeroth left derived functor $L_0^I$ of $I$-adic completion as the 'correct' notion.

**Definition 5.1.**
- We say that a module $M$ is $L_0^I$-*complete* if the natural map $M \to L_0^I M$ is an isomorphism.
- We say that a dg-module $N$ is $L_0^I$-*complete* if its underlying graded module is $L_0^I$-complete.

We write $\mathrm{Mod}_R$ for the category of dg-$R$-modules, and $\mathrm{Mod}_R^\wedge$ for the full subcategory of $L_0^I$-complete dg-modules. We denote the internal hom of $R$-modules by $\underline{\mathrm{Hom}}_R(-,-)$.

**Lemma 5.2.**
(a) *The category $\mathrm{Mod}_R^\wedge$ is abelian, and the inclusion functor $i\colon \mathrm{Mod}_R^\wedge \to \mathrm{Mod}_R$ is exact. In particular, the homology of an $L_0^I$-complete dg-module is $L_0^I$-complete.*
(b) *The inclusion functor is right adjoint to the L-completion functor $L_0^I$.*
(c) *The category $\mathrm{Mod}_R^\wedge$ has all limits and colimits.*

*Proof.* The proofs of (a) and (b) can be found in [31, A.6(e), A.6(f)]. Their proofs depend only upon the fact that $L_0^I$ is right exact and the existence of a long exact sequence of derived functors. Therefore, the restriction to local rings and regular ideals made in [31] does not affect the stated results. It follows from (b) that limits of $L_0^I$-complete modules are calculated in $\mathrm{Mod}_R$, and that colimits of $L_0^I$-complete modules are calculated by $L_0^I$-completion of the colimit in $\mathrm{Mod}_R$. □

**Proposition 5.3.**
(a) *If $N$ is $L_0^I$-complete, then $\underline{\mathrm{Hom}}_R(M,N)$ is $L_0^I$-complete.*
(b) *The category $\mathrm{Mod}_R^\wedge$ is closed symmetric monoidal with tensor product $L_0^I(M \otimes_R N)$ and internal hom $\underline{\mathrm{Hom}}_R(M,N)$.*

*Proof.* By taking a free presentation $R^{J_1} \to R^{J_0} \to M \to 0$, we obtain an exact sequence
$$0 \to \underline{\mathrm{Hom}}_R(M,N) \to \prod_{J_0} N \to \prod_{J_1} N$$
which proves (a), since $L_0^I$-complete modules are closed under products and kernels.

For (b) we follow the argument of Rezk [38, 6.2]. We first prove that the map $L_0^I(M \otimes_R N) \to L_0^I(L_0^I M \otimes_R N)$ induced by $\eta_M \colon M \to L_0^I M$ is an isomorphism. It is enough to check that for any $L_0^I$-complete module $C$, the map
$$\underline{\mathrm{Hom}}_R(L_0^I(L_0^I M \otimes_R N), C) \to \underline{\mathrm{Hom}}_R(L_0^I(M \otimes_R N), C)$$
is an isomorphism. By adjunction, it is an isomorphism if and only if the induced map
$$\underline{\mathrm{Hom}}_R(L_0^I M, \underline{\mathrm{Hom}}_R(N,C)) \to \underline{\mathrm{Hom}}_R(M, \underline{\mathrm{Hom}}_R(N,C))$$
is. This now follows as $\underline{\mathrm{Hom}}_R(N,C)$ is $L_0^I$-complete by part (a). Therefore $L_0^I(M \otimes_R N) \to L_0^I(L_0^I M \otimes_R N)$ is an isomorphism. By symmetry, we also have that $L_0^I(M \otimes_R N) \to L_0^I(M \otimes_R L_0^I N)$ is an isomorphism, and therefore so is $L_0^I(M \otimes_R N) \to L_0^I(L_0^I M \otimes_R L_0^I N)$. This completes the proof of (b). □

We will also be concerned with a homotopical version of completion that we shall now recall. For any $x \in R$, we define the *unstable Koszul complex*

$$K(x) = \mathrm{fib}(\Sigma^{|x|}R \xrightarrow{\cdot x} R),$$

and the *stable Koszul complex*

$$K_\infty(x) = \mathrm{fib}(R \to R[1/x])$$

where the fibre is taken in the category of dg-modules. For an ideal $I = (x_1, \ldots, x_n)$ we put

$$K(I) = K(x_1) \otimes_R \cdots \otimes_R K(x_n) \quad \text{and} \quad K_\infty(I) = K_\infty(x_1) \otimes_R \cdots \otimes_R K_\infty(x_n).$$

If no confusion is likely to arise, we suppress notation for the ideal and write $K$ for the unstable Koszul complex and $K_\infty$ for the stable Koszul complex. We will also write $\mathrm{Hom}_R(-,-)$ for the derived internal hom functor. We say that a dg-module $M$ is *derived complete* if the natural map $M \to \mathrm{Hom}_R(K_\infty, M) =: \Lambda_I M$ is a quasi-isomorphism. Then the *nth local homology* of $M$ is defined to be $H_n^I(M) = H_n(\Lambda_I M)$.

**Definition 5.4.** Let $I = (x_1, \ldots, x_n)$ be a finitely generated homogeneous ideal. For all $s \in \mathbb{N}$ and $x \in R$, we put

$$K_s(x) = \mathrm{fib}(\Sigma^{s|x|}R \xrightarrow{\cdot x^s} R) \quad \text{and} \quad K_s(I) = K_s(x_1) \otimes_R \cdots \otimes_R K_s(x_n).$$

We say that $I$ is generated by the *weakly pro-regular* sequence $(x_1 \ldots, x_n)$ if the inverse system $(H_k(K_s(I)))_s$ is pro-zero for all $k \neq 0$. That is, for each $s \in \mathbb{N}$ there is $m \geq s$ such that the natural map

$$H_k(K_m(I)) \to H_k(K_s(I))$$

is zero.

Note that if $R$ is Noetherian then any finitely generated ideal is weakly pro-regular [37, 4.34]. Indeed this is true even when weakly pro-regular is replaced by pro-regular [21].

**Theorem 5.5.** *Let $R$ be a graded commutative ring and let $I$ be a finitely generated homogeneous ideal that is generated by a weakly pro-regular sequence. Then for all dg-modules $M$, there is a natural quasi-isomorphism*

$$\mathrm{tel}^L_{I,M} \colon \Lambda_I(M) \xrightarrow{\sim} \mathbb{L}(-)^\wedge_I(M)$$

*between the derived completion functor and the total left derived functor of $I$-adic completion, making the diagram*

$$\begin{array}{ccc}
& M & \\
\swarrow & & \searrow \\
\Lambda_I M & \xrightarrow{\mathrm{tel}^L_{I,M}} & \mathbb{L}(-)^\wedge_I(M)
\end{array}$$

*commute. Moreover, taking homology on both sides we get*

$$H_*^I(M) \cong L_*^I(M).$$

*Proof.* Greenlees-May proved that if $R$ has bounded torsion and $I$ is pro-regular then $H_*^I M \cong L_*^I M$, see [21, 2.5]. Schenzel [39, 1.1] proved the above result for ideals generated by weakly pro-regular sequences and bounded complexes with $R$ bounded torsion. Finally, Porta-Shaul-Yekutieli [37, 5.25] removed the hypothesis that $R$ has bounded torsion and extended the result to unbounded complexes. $\square$

As an application, we prove the following result which we will use in the construction of an Adams spectral sequence, see Theorem 10.6.

**Proposition 5.6.** *Let $R$ be a graded commutative ring and let $I$ be a finitely generated homogeneous ideal that is generated by a weakly pro-regular sequence.*

(a) *If $M$ is an $L_0^I$-complete module and $P_\bullet \to M$ is a projective resolution of $M$, then $L_0^I P_\bullet \to M$ is a projective resolution in $L_0^I$-complete modules.*

(b) *Write $\widehat{\mathrm{Ext}}_R$ for the Ext-groups in the abelian category of $L_0^I$-complete modules. Then*
$$\widehat{\mathrm{Ext}}_R(M,N) \cong \mathrm{Ext}_R(M,N)$$
*for all $L_0^I$-complete modules $M$ and $N$.*

*Proof.* Given an $L_0^I$-complete module $M$, choose a projective resolution $P_\bullet \to M$ in $R$-modules. Since $L_0^I$ is left adjoint to the inclusion, $L_0^I P_\bullet$ is a complex of projective $L_0^I$-complete modules. We now show that $L_0^I P_\bullet \to M$ is a projective resolution in $L_0^I$-complete modules. Note that $\Lambda_I P_\bullet \to \Lambda_I M$ is a quasi-isomorphism. Using Theorem 5.5 and [21, 4.1] we have that $L_0^I P_\bullet \simeq \Lambda_I P_\bullet$ and $M \simeq \Lambda_I M$. Therefore $L_0^I P_\bullet \to M$ is a projective resolution in $L_0^I$-complete modules. By adjunction we deduce that for all $L_0^I$-complete modules $M$ and $N$ we have
$$\widehat{\mathrm{Ext}}_R(M,N) = H_*(\mathrm{Hom}_R(L_0^I P_\bullet, N)) \cong H_*(\mathrm{Hom}_R(P_\bullet, N)) = \mathrm{Ext}_R(M,N).$$
□

**Corollary 5.7.** *Let $R$ be a graded commutative ring and let $I$ be a finitely generated homogeneous ideal that is generated by a weakly pro-regular sequence. Let $\{M_i\}$ be a collection of $L_0^I$-complete modules. Then $L_n^I(\bigoplus M_i) = 0$ for all $n \geq 1$*[3].

*Proof.* Let $P_i$ be a projective resolution of $M_i$ in $R$-modules. We have that
$$L_n^I(\oplus M_i) = H_n(L_0^I(\oplus P_i)) \cong H_n(\oplus P_i)$$
where the isomorphism follows from Proposition 5.6. Since direct sums are exact, this is isomorphic to $\bigoplus L_n^I M_i$ which is zero for all $n \geq 1$ by [21, 4.1]. □

## 6. An abelian model for derived completion

In this section we use the language of model categories to show that the category of $L_0^I$-complete modules forms an abelian model for derived complete modules, see Theorem 6.10. Our result can be thought as "dual" to the fact that $I$-power torsion modules forms an abelian model for derived torsion modules [25, §5]. We will be working under the following:

**Hypothesis 6.1.** We will assume our ideal $I$ to be generated by a weakly pro-regular sequence and we continue to write $K$ for its associated unstable Koszul complex.

In order to prove our main result we need to consider model structures on the categories of interest. Recall that the category of dg-modules $\mathrm{Mod}_R$ has a projective model structure in which the weak equivalences are the quasi-isomorphisms, the fibrations are the epimorphisms and the cofibrations are the monomorphisms which have dg-projective cokernel and are split on the underlying graded modules, see [13, 3.3] and [1, 3.15]. A dg-module $M$ is said to be *dg-projective* if $\mathrm{Hom}_R(P, -)$ preserves surjective quasi-isomorphisms. It is important to note that any dg-projective module is (graded) projective, but the converse is not generally true, see [3, 9.6.1].

**Lemma 6.2.** *If $P$ is dg-projective, then $L_n^I P = 0$ for all $n \geq 1$. Moreover, there is a natural quasi-isomorphism $\Lambda_I P \xrightarrow{\sim} L_0^I P$.*

*Proof.* This is the trivial case of Theorem 5.5. □

We will now put a projective model structure on $L_0^I$-complete modules following Rezk's unpublished note [38, 10.2].

**Lemma 6.3.**

(a) *The functor $L_0^I$ takes cofibrations in $\mathrm{Mod}_R$ to morphisms which have the left lifting property with respect to surjective quasi-isomorphisms of $L_0^I$-complete modules.*

---
[3]This result fails if the ideal is generated by a sequence of more than length one; see the corrigendum at the end of this paper

(b) *The functor $L_0^I$ takes acyclic cofibrations in $\mathrm{Mod}_R$ to morphisms which have the left lifting property with respect to surjections of $L_0^I$-complete modules.*

(c) *If $M \to N$ is a cofibration in $\mathrm{Mod}_R$, the homology $H_*N$ is $L_0^I$-complete and $M \to L_0^I M$ is a quasi-isomorphism, then $N \to L_0^I N$ is a quasi-isomorphism.*

*Proof.* Part (a) and (b) follow from the lifting properties in $\mathrm{Mod}_R$. For part (c), note that by definition $M \to N$ is an injection with dg-projective cokernel $P$ so we have a diagram

$$
\begin{array}{ccccccccc}
0 & \longrightarrow & M & \longrightarrow & N & \longrightarrow & P & \longrightarrow & 0 \\
& & \sim\downarrow & & \downarrow & & \downarrow & & \\
0 & \longrightarrow & L_0^I M & \longrightarrow & L_0^I N & \longrightarrow & L_0^I P & \longrightarrow & 0
\end{array}
$$

in which the top row is exact. By Lemma 6.2 we have that $L_i^I P = 0$ for $i \geq 1$ so the long exact sequence of derived functors collapses to a short exact sequence. Therefore, the bottom row is exact too. Since $L_0^I M$ is $L_0^I$-complete, the homology $H_*M \cong H_*L_0^I M$ is $L_0^I$-complete by Lemma 5.2(a), and so $H_*P$ is $L_0^I$-complete too. Now consider the spectral sequence [22, 3.3]

$$E_{p,q}^2 = (L_p^I(H_*P))_q \implies H_{p+q}(\Lambda_I P).$$

If the homology groups are $L_0^I$-complete, then the spectral sequence collapses by [21, 4.1] to give a quasi-isomorphism $P \to \Lambda_I P$. Therefore $P \to L_0^I P$ is a quasi-isomorphism by Lemma 6.2. Hence $N \to L_0^I N$ is a quasi-isomorphism as required. □

**Proposition 6.4.** *There is a model structure on $\mathrm{Mod}_R^\wedge$ in which the weak equivalences are the quasi-isomorphisms, the fibrations are the surjections, and the cofibrations are the maps with the left lifting property with respect to the acyclic fibrations. Furthermore, the adjunction*

$$L_0^I : \mathrm{Mod}_R \rightleftarrows \mathrm{Mod}_R^\wedge : i$$

*is Quillen.*

*Proof.* The only parts that need elaboration are the factorization axiom and the lifting axiom. Firstly we prove the factorization axiom.

Let $f \colon M \to N$ in $\mathrm{Mod}_R^\wedge$. Take a factorization $M \xrightarrow{i} D \xrightarrow{p} N$ in $\mathrm{Mod}_R$ where one of $i$ or $p$ is acyclic. Since $L_0^I$ is left adjoint to the inclusion, maps $L_0^I D \to N$ are in bijection with maps $D \to N$. Therefore, there is a unique $q \colon L_0^I D \to N$ making the square

$$
\begin{array}{ccc}
M & \longrightarrow & L_0^I D \\
i\downarrow & \nearrow \alpha & \downarrow q \\
D & \xrightarrow{p} & N
\end{array}
$$

commute. Note that $q$ is a fibration since $q \cong L_0^I p$ and $L_0^I$ preserves surjections.

If $p$ is acyclic, Lemma 6.3(c) shows that $\alpha$ is a quasi-isomorphism since $H_*D \cong H_*N$, and hence by the two-of-three property, $q$ is a weak equivalence. Lemma 6.3(a) shows that the factorization $f = q(\alpha i)$ is a factorization into a map with the left lifting property with respect to acyclic fibrations, followed by an acylic fibration. This completes the first part of the proof of the factorization axiom.

For the other part we suppose that $i$ is a weak equivalence. Since $\alpha i \cong L_0^I(i)$, Lemma 6.3(b) shows that $\alpha i$ has the left lifting property with respect to fibrations in $\mathrm{Mod}_R^\wedge$. Lemma 6.3(c) shows that $\alpha$ is a quasi-isomorphism since $H_*D \cong H_*M$. Therefore $f = q(\alpha i)$ is a factorization into a weak equivalence with the left lifting property with respect to fibrations followed by an fibration, which completes the proof of the factorization axiom.

For the lifting axiom, we note that one part is by definition. For the other part, we use the standard method of the retract argument. Consider the square

$$\begin{array}{ccc} A & \longrightarrow & X \\ \downarrow i & & \downarrow f \\ B & \longrightarrow & Y \end{array}$$

in which $i$ is an acyclic cofibration and $f$ is a fibration. Factor $i$ into a map with the left lifting property with respect to fibrations followed by a fibration to give $A \xrightarrow{j} C \xrightarrow{p} B$. Since $j$ has the left lifting property with respect to fibrations, there is a lift $g \colon C \to X$.

As $i$ and $j$ are weak equivalences, $p$ is an acyclic fibration. Since $i$ has the left lifting property with respect to acyclic fibrations, there exists a lift $h \colon B \to C$. Therefore $gh \colon B \to X$ gives the required lift in the square.

It is clear that the adjunction is Quillen by the definition of the weak equivalences and fibrations. $\square$

**Remark 6.5.** One might first think of attempting to prove the existence of this model structure via right inducing it from $\mathrm{Mod}_R$. However, in order to be able to do this, we need to know that the inclusion $i \colon \mathrm{Mod}_R^\wedge \to \mathrm{Mod}_R$ preserves filtered colimits. This is false; take $R = \mathbb{Z}$ and $I = (p)$ and consider the direct system $\mathbb{Z}_p \xrightarrow{p} \mathbb{Z}_p \xrightarrow{p} \ldots$. Then the colimit in the category of abelian groups is $\mathbb{Q}_p$, while the colimit in the category of $L_0^I$-complete abelian groups is $L_0^{(p)}(\mathbb{Q}_p)$ which is zero.

**Proposition 6.6.** *The model structure on* $\mathrm{Mod}_R^\wedge$ *is symmetric monoidal.*

*Proof.* The category of $L_0^I$-complete modules is closed symmetric monoidal with tensor product given by $L_0^I(M \otimes N)$; see Proposition 5.3.

Let $M \to N$ and $X \to Y$ be fibrations in $\mathrm{Mod}_R^\wedge$. Since the inclusion $i \colon \mathrm{Mod}_R^\wedge \to \mathrm{Mod}_R$ preserves limits, we have that the pullback product map is

$$\underline{\mathrm{Hom}}_R(iN, iX) \to \underline{\mathrm{Hom}}_R(iM, iX) \times_{\underline{\mathrm{Hom}}_R(iM, iY)} \underline{\mathrm{Hom}}_R(iN, iY).$$

Since $\mathrm{Mod}_R$ is a symmetric monoidal model category and $i$ is right Quillen, the pullback product map is a fibration. A similar proof shows that the pullback product of a fibration with an acyclic fibration is an acyclic fibration. The unit axiom is immediate since the unit in $\mathrm{Mod}_R^\wedge$ is $L_0^I R$ which is cofibrant as $R$ is cofibrant in $\mathrm{Mod}_R$. $\square$

We need a model category modelling the homotopy theory of derived complete dg-modules. The left Bousfield localization of $R$-modules at the unstable Koszul complex is such a model category by the following result.

**Lemma 6.7** ([22, 4.2]). *There is an equivalence of categories*

$$\mathrm{h}L_K \mathrm{Mod}_R \simeq \Lambda \mathrm{Mod}_R$$

*where* $\Lambda \mathrm{Mod}_R$ *denotes the full subcategory of the derived category of dg-modules consisting of derived complete dg-modules.*

We relate the model category of $L_0^I$-complete modules to derived complete modules. We will use these results to show that cofree $G$-spectra have an abelian model in terms of $L_0^I$-complete modules.

**Lemma 6.8.** *There is a symmetric monoidal Quillen adjunction*

$$L_0^I \colon L_K(\mathrm{Mod}_R) \rightleftarrows \mathrm{Mod}_R^\wedge \colon i.$$

*Proof.* The cofibrations in $L_K \mathrm{Mod}_R$ are the same as the cofibrations in $\mathrm{Mod}_R$ so they are preserved since $L_0^I \colon \mathrm{Mod}_R \to \mathrm{Mod}_R^\wedge$ is left Quillen. Now suppose that $f \colon M \to N$ is an acyclic cofibration in $L_K(\mathrm{Mod}_R)$ so that the cokernel $C$ is dg-projective. In particular, $K \otimes C$ and $\mathrm{Hom}_R(K, C)$ are acyclic as $K$ is self-dual up to suspension. We also know that $K_\infty$ is built from $K$ so $\Lambda_I C = \mathrm{Hom}_R(K_\infty, C)$ is acyclic as well. By Lemma 6.2, we have $\Lambda_I C \simeq L_0^I C$ and so $L_0^I M \to L_0^I N$ is a quasi-isomorphism. This is a symmetric monoidal Quillen adjunction since $L_0^I$ is strong monoidal by Lemma 5.2, and the unit in $L_K(\mathrm{Mod}_R)$ is cofibrant. $\square$

Before we can prove that the above Quillen adjunction is actually a Quillen equivalence, we need the following:

**Lemma 6.9.** *For any dg-module* $M$, *the natural map* $K \otimes M \to \Lambda_I(K \otimes M)$ *is a quasi-isomorphism.*

*Proof.* There is a fibre sequence $K_\infty \to R \to \check{C}R$ where $\check{C}R = \Sigma \ker(K_\infty \to R)$ is the Čech complex. This gives rise to another fibre sequence
$$\mathrm{Hom}_R(K_\infty, N) \leftarrow N \leftarrow \mathrm{Hom}_R(\check{C}R, N)$$
for any dg-module $N$. Now let $I = (x_1, \ldots, x_n)$. Note that $\check{C}R$ is finitely built from $R[\frac{1}{x_i}]$ and that the multiplication map $x_i \colon K \to K$ is null-homotopic. Thus $\mathrm{Hom}_R(\check{C}R, K \otimes M) \simeq 0$ and $K \otimes M$ is derived complete. □

We can now prove that $L_0^I$-complete modules are a model for derived complete modules.

**Theorem 6.10.** *There is a symmetric monoidal Quillen equivalence*
$$L_0^I : L_K(\mathrm{Mod}_R) \rightleftarrows \mathrm{Mod}_R^\wedge : i.$$

*Proof.* We now show that this Quillen adjunction is in fact a Quillen equivalence. Let $P$ be cofibrant (i.e., dg-projective) in $L_K(\mathrm{Mod}_R)$ and $M$ be fibrant in the category of $L_0^I$-complete $R$-modules. We must show that $L_0^I P \to M$ is a quasi-isomorphism if and only if $K \otimes P \to K \otimes M$ is a quasi-isomorphism.

Firstly, if $L_0^I P \to M$ is a quasi-isomorphism, then $K \otimes L_0^I P \to K \otimes M$ is a quasi-isomorphism since $K$ is homotopically flat. Now note that there is a weak equivalence $K \otimes \Lambda_I P \xrightarrow{\sim} \Lambda_I(K \otimes P)$ since $K$ is small. By Lemma 6.2, $K \otimes L_0^I P \simeq \Lambda_I(K \otimes P)$ as $P$ is projective. Hence $K \otimes L_0^I P \simeq K \otimes P$ by Lemma 6.9. We conclude that $K \otimes P \to K \otimes M$ is a quasi-isomorphism as required.

Conversely, if $K \otimes P \to K \otimes M$ is a quasi-isomorphism then $\mathrm{Hom}_R(K, P) \to \mathrm{Hom}_R(K, M)$ is too since $K$ is self-dual up to suspension. Since $K_\infty$ is built from $K$, we also deduce $\mathrm{Hom}_R(K_\infty, P) \to \mathrm{Hom}_R(K_\infty, M)$ is a quasi-isomorphism. It follows that $\Lambda_I P \to \Lambda_I M$ is a quasi-isomorphism. By Lemma 6.2, we have $L_0^I P \simeq \Lambda_I P$ and $M \simeq \Lambda_I M$. Hence $L_0^I P \to M$ is a quasi-isomorphism. □

As a consequence we obtain the following corollary which extends [17, 6.15] to non-Noetherian rings.

**Corollary 6.11.** *A dg-module $M$ is derived complete if and only if its homology $H_*M$ is $L_0^I$-complete.*

*Proof.* Let $M$ be derived complete. By Theorem 6.10, $M$ is quasi-isomorphic to its $L_0^I$-completion $L_0^I M$. As the homology of an $L_0^I$-complete object is still $L_0^I$-complete by Lemma 5.2, we deduce that $M$ has $L_0^I$-complete homology. Conversely, suppose that $M$ is a module with $L_0^I$-complete homology. The spectral sequence [22, 3.3]
$$E^2_{p,q} = (L_p^I H_* M)_q \implies H_{p+q}(\Lambda_I M)$$
collapses by [21, 4.1], showing that $M \to \Lambda_I M$ is a quasi-isomorphism. □

## 7. The category of rational cofree $G$-spectra

From now on we will be working rationally. This means that all spectra are rationalized without comment and all homology and cohomology theories will be unreduced and with rational coefficients.

**Notation 7.1.** Fix $G$ a compact Lie group. We denote by $\mathrm{Sp}_G$ the model category of rational orthogonal $G$-spectra with the rational $G$-stable model structure, which is a compactly generated, stable, symmetric monoidal model category, see [34, III.7.6]. We write $\wedge$ for the tensor product and $F(-,-)$ for the internal hom functor. We also write $\mathrm{hSp}_G$ for its associated homotopy category.

**Definition 7.2.** A $G$-spectrum $X$ is said to be *cofree* if the natural map $X \to F(EG_+, X)$ is an isomorphism in the homotopy category. We denote by $\mathrm{hSp}_G^{\mathrm{cofree}}$ the full subcategory of $\mathrm{hSp}_G$ of cofree $G$-spectra.

**Lemma 7.3.** *There is a natural equivalence*
$$\mathrm{h}L_{EG_+}\mathrm{Sp}_G \simeq \mathrm{hSp}_G^{\mathrm{cofree}}.$$
*Furthermore, $L_{EG_+}\mathrm{Sp}_G$ is a symmetric monoidal model category.*

*Proof.* A fibrant replacement functor in $L_{EG_+}\mathrm{Sp}_G$ is given by $F(EG_+, R(-))$ where $R$ is the fibrant replacement in $\mathrm{Sp}_G$. Therefore, the collection of bifibrant objects in $L_{EG_+}\mathrm{Sp}_G$ is equivalent to the full subcategory of cofree $G$-spectra. The model category $L_{EG_+}\mathrm{Sp}_G$ is symmetric monoidal by Proposition 3.11. □

## 8. The symmetric monoidal equivalence: connected case

In this section we fix a connected compact Lie group $G$. We aim to find an algebraic model for the category of rational cofree $G$-spectra. There are several steps needed. Recall that our model for cofree $G$-spectra is the homological localization $L_{EG_+}\mathrm{Sp}_G$.

**Step 1.** Any cofree $G$-spectrum is naturally a module over $F(EG_+, EG_+)$ which is equivalent to the complex orientable commutative ring $G$-spectrum $DEG_+ = F(EG_+, S^0)$. Restriction and extension of scalars along the unit map $S^0 \to DEG_+$ induces a symmetric monoidal Quillen adjunction

$$DEG_+ \wedge - : L_{EG_+}(\mathrm{Sp}_G) \rightleftarrows L_{EG_+}(\mathrm{Mod}_{DEG_+}) : U$$

between the localizations, since $DEG_+ \wedge EG_+ \simeq EG_+$. By the Left Localization Principle this is a symmetric monoidal Quillen equivalence, since the unit is an $EG_+$-equivalence and $U$ preserves non-equivariant equivalences.

**Remark 8.1.** This is a special case of Proposition 4.4 and Example 4.7.

**Step 2.** We can now take categorical fixed points to remove equivariance. As a functor from $G$-spectra to non-equivariant spectra, the categorical fixed points is right adjoint to the inflation functor. Using [42, §3.3] we have a symmetric monoidal Quillen adjunction

$$(-)^G : \mathrm{Mod}_{DEG_+} \leftrightarrows \mathrm{Mod}_{DBG_+} : DEG_+ \otimes_{DBG_+} -$$

between the categories of modules. Note that we suppress notation for the inflation functor. A more detailed discussion of this adjunction can be found in [27].

Since $G$ is connected, $DEG_+$ generates $\mathrm{Mod}_{DEG_+}$ by [20, 3.1] and so the counit is an equivalence on all objects as it is an equivalence on $DEG_+$ and the fixed points functor preserves sums. By [20, 3.3], the fixed points functor sends non-equivariant equivalences to $BG_+$-equivalences, so the Left Localization Principle applies and we get a symmetric monoidal Quillen equivalence

$$(-)^G : L_{EG_+}\mathrm{Mod}_{DEG_+} \leftrightarrows L_{BG_+}\mathrm{Mod}_{DBG_+} : DEG_+ \otimes_{DBG_+} -.$$

**Step 3.** We now apply Shipley's theorem [45, 2.15] (see also [47, 7.2]) which gives a symmetric monoidal Quillen equivalence

$$\Theta \colon \mathrm{Mod}_{DBG_+} \simeq_Q \mathrm{Mod}_{\Theta DBG_+}$$

where $\Theta DBG_+$ is a commutative dga with the property that $H_*(\Theta DBG_+) = \pi_*(DBG_+) = H^*BG$. It follows that there is a symmetric monoidal Quillen equivalence

$$L_{BG_+}\mathrm{Mod}_{DBG_+} \simeq_Q L_{\Theta BG_+}\mathrm{Mod}_{\Theta DBG_+}$$

where $H_*(\Theta BG_+) \cong \pi_*(BG_+) \cong H_*BG$.

**Step 4.** Since $H^*BG$ is a polynomial ring it is strongly intrinsically formal as a commutative dga. In other words, for any commutative dga $R$ with $H_*R \cong H^*BG$, there is a quasi-isomorphism $H^*BG \to R$. Therefore, taking cycle representatives we have a quasi-isomorphism $z \colon H^*BG \to \Theta DBG_+$. We also need the following result to identify $\Theta BG_+$.

**Lemma 8.2.** *There is a natural weak equivalence $\Theta BG_+ \to H_*BG$.*

*Proof.* Write $(-)^\vee = \mathrm{Hom}_\mathbb{Q}(-, \mathbb{Q})$ and note that it is exact. There is a canonical map $\Theta BG_+ \to (\Theta BG_+)^{\vee\vee}$ which is a quasi-isomorphism since the homotopy groups of $BG_+$ are degreewise finite. There is a natural map $\Theta DBG_+ \to (\Theta BG_+)^\vee$ obtained as the transpose of the natural composite

$$\Theta BG_+ \otimes \Theta DBG_+ \to \Theta(BG_+ \wedge DBG_+) \to \mathbb{Q}.$$

Since $\Theta$ gives a symmetric monoidal equivalence of homotopy categories, the natural map $\Theta DBG_+ \to (\Theta BG_+)^\vee$ is a weak equivalence.

Since $DBG_+$ is a commutative $H\mathbb{Q}$-algebra, $\Theta DBG_+$ is a commutative dga by [45, 1.2]. As $H^*BG$ is strongly intrinsically formal as a commutative dga, there exists a quasi-isomorphism $H^*BG \to \Theta DBG_+$. Putting all this together, we have quasi-isomorphisms
$$\Theta BG_+ \to (\Theta BG_+)^{\vee\vee} \to (H^*BG)^\vee \to H_*BG.$$
$\square$

Extension and restriction of scalars along the map $z\colon H^*BG \to \Theta DBG_+$
$$\mathrm{Mod}_{\Theta DBG_+} \xleftarrow[z^*]{\Theta DBG_+ \otimes_{H^*BG} -} \mathrm{Mod}_{H^*BG}$$
is a symmetric monoidal Quillen equivalence since chain complexes satisfies Quillen invariance of modules. Therefore we have a symmetric monoidal Quillen equivalence
$$L_{H_*BG}\mathrm{Mod}_{\Theta DBG_+} \simeq_Q L_{H_*BG}\mathrm{Mod}_{H^*BG}.$$

**Step 5.** It remains to internalize the localization. Let $I$ be the augmentation ideal of $H^*BG$ and let $K$ denote its unstable Koszul complex.

**Proposition 8.3.** *The homology $H_*BG$ finitely builds $K$ and $K$ builds $H_*BG$.*

*Proof.* Suppose that $H^*BG = \mathbb{Q}[x_1,...,x_n]$. There is a cofibre sequence
$$\Sigma^{|x_1|}\mathbb{Q}[x_1,...,x_n] \xrightarrow{\cdot x_1} \mathbb{Q}[x_1,...,x_n] \to \Sigma K(x_1)$$
and applying $\mathrm{Hom}_\mathbb{Q}(-,\mathbb{Q})$ gives the cofibre sequence
$$H_*BG \to \Sigma^{-|x_1|}H_*BG \to \Sigma K(x_1)^\vee.$$
Since $K(x_1)$ is self-dual up to suspension, this shows that $K(x_1)$ is finitely built from $H_*BG$. A repeated argument using the cofibre sequence $\Sigma^{|x_i|}K_{i-1} \to K_{i-1} \to K_i$ where $K_i = K(x_1,...,x_i)$ and $K_0 = H^*BG$ shows that $K$ is finitely built from $H_*BG$.

Conversely, since $H_*BG$ is torsion it is built by $K$ as $K$ generates torsion modules [24, 8.7]. $\square$

Therefore, a map is a $H_*BG$-equivalence if and only if it is a $K$-equivalence. It follows that
$$L_{H_*BG}\mathrm{Mod}_{H^*BG} = L_K\mathrm{Mod}_{H^*BG}.$$

Combining all the statements of this section with Theorem 6.10 gives the following result.

**Theorem 8.4.** *Let $G$ be a connected compact Lie group and $I$ be the augmentation ideal of $H^*BG$. Then there is a symmetric monoidal Quillen equivalence*
$$L_{EG_+}\mathrm{Sp}_G \simeq_Q \mathrm{Mod}^\wedge_{H^*BG}$$
*between rational cofree $G$-spectra and $L_0^I$-complete dg-$H^*BG$-modules. In particular, there is a tensor-triangulated equivalence*
$$\text{cofree } G\text{-spectra} \simeq_\triangle \mathcal{D}(L_0^I\text{-complete } H^*BG\text{-modules}).$$

9. The symmetric monoidal equivalence: non-connected case

In this section we extend the algebraic model for cofree $G$-spectra from connected $G$ to any compact Lie group. The blueprint is the same as for the connected case, however some extra care is required which arises from taking categorical fixed points. We fix a compact Lie group $G$ with identity component $N$ and component group $W = G/N$, and write $r$ for the rank of $G$.

9.1. **Skewed Model Categories.** We recall some results about model categories with a action of a finite group $W$ from [32, §5.2] and [5, §7]. For any cofibrantly generated model category $\mathcal{C}$, we denote by $\mathcal{C}[W] = \text{Fun}(BW, \mathcal{C})$ the category of objects of $\mathcal{C}$ with a $W$-action. We endow $\mathcal{C}[W]$ with the projective model structure where the weak equivalence and fibrations are created by the forgetful functor $\mathcal{C}[W] \to \mathcal{C}$. We will need the following result:

**Lemma 9.1** ([32, 5.3]). *There is a symmetric monoidal Quillen equivalence $L_{EW_+}\text{Sp}_W \simeq_Q \text{Sp}[W]$.*

More generally, we can consider the category $\mathbb{E}W$ with objects the elements of $W$ and a unique morphism connecting each pair of objects. Let $\mathcal{C}$ be a category with a $W$-action, that is, with functors $w_* \colon \mathcal{C} \to \mathcal{C}$ for each $w \in W$ satisfying $(ww')_* = w_* w'_*$ and $e_* = 1$. The category of objects of $\mathcal{C}$ with a *skewed $W$-action* is the category of equivariant functors $\mathbb{E}W \to \mathcal{C}$ and equivariant natural transformations, which we denote by $\mathcal{C}[\widetilde{W}]$. Note that if the $W$-action on $\mathcal{C}$ is trivial, then $\mathcal{C}[\widetilde{W}]$ is equivalent to $\mathcal{C}[W]$. We say that an adjunction between categories with a $W$-action is a *$W$-adjunction* if both the functors are $W$-equivariant and the unit and counit are $W$-equivariant natural transformations. We say that a model category $\mathcal{C}$ with a $W$-action is *skewable* if $w_* \colon \mathcal{C} \to \mathcal{C}$ is left Quillen for each $w \in W$. Note that $w_* \colon \mathcal{C} \to \mathcal{C}$ is left adjoint to $w_*^{-1}$, so equivalently, we could ask for $w_*$ to be right Quillen for all $w \in W$.

**Lemma 9.2.**

(a) *Let $\mathcal{C}$ be a skewable, symmetric monoidal, cofibrantly generated model category with a $W$-action. Then $\mathcal{C}[\widetilde{W}]$ admits a closed symmetric monoidal structure and a projective model structure making it into a symmetric monoidal model category.*

(b) *Let $\mathcal{C}$ and $\mathcal{D}$ be skewable, symmetric monoidal model categories. Suppose that $\mathcal{C} \rightleftarrows \mathcal{D}$ is a $W$-adjunction which is a symmetric monoidal Quillen equivalence. Then we have a symmetric monoidal Quillen equivalence*
$$\mathcal{C}[\widetilde{W}] \simeq_Q \mathcal{D}[\widetilde{W}].$$

*Proof.* One can check that $\mathcal{C}[\widetilde{W}]$ is a symmetric monoidal model category in which the weak equivalences and fibrations are determined levelwise, and that Quillen equivalences extend to the skewed model category; see [5, §7.3] for the case $W = C_2$. $\square$

9.2. **The algebraic model.** The component group $W$ acts on $N$ by conjugation and hence on its cohomology $H^*BN$. We write $H^*\widetilde{BN}$ for the polynomial ring $H^*BN$ equipped with this $W$-action. Accordingly, the model category $\text{Mod}_{H^*\widetilde{BN}}$ inherits a $W$-action as follows. For $w \in W$ and a $H^*\widetilde{BN}$-module $M$, we define $w_*M$ to be the same underlying abelian group as $M$ but with module structure now defined by $r \cdot m := (wr)m$ for $r \in H^*\widetilde{BN}$ and $m \in M$. This model category is skewable since the action preserves weak equivalences and fibrations. Therefore, we can consider the model category $\text{Mod}_{H^*\widetilde{BN}}[\widetilde{W}]$ of modules with a skewed $W$-action. More explicitly, we can identify this category with the category of modules over the skewed ring $H^*\widetilde{BN}[W]$, that is, the ring whose elements are formal linear sums $\sum_{w \in W} x_w w$ where $x_w \in H^*\widetilde{BN}$, with pointwise addition and multiplication defined by
$$(xw) \cdot (x'w') = (x(w \cdot x'))(ww') \quad \text{for } w, w' \in W \text{ and } x, x' \in H^*\widetilde{BN}.$$

We now turn to define a suitable notion of $L_0^I$-completion for a module over the skewed ring.

**Definition 9.3.** Let $I$ denote the augmentation ideal of $H^*BN$. We say that a dg-$H^*\widetilde{BN}[W]$-module $M$ is $L_0^I$-*complete* if $M$ is $L_0^I$-complete as a $H^*BN$-module. We denote by $\text{Mod}^\wedge_{H^*\widetilde{BN}[W]}$ the category of $L_0^I$-complete dg modules over the skewed ring.

**Lemma 9.4.**

(a) *The category of left $H^*\widetilde{BN}[W]$-modules admits a closed symmetric monoidal structure and a projective model structure making it into a symmetric monoidal model category.*

(b) *The category of $L_0^I$-complete left $H^*\widetilde{BN}[W]$-modules is abelian and is a symmetric monoidal model category with the projective model structure.*

*Proof.* The results follow from the previous sections and Lemma 9.2 by noticing that the category of ($L_0^I$-complete) $H^*\widetilde{BN}[W]$-modules is equivalent to $\mathcal{C}[\widetilde{W}]$ for $\mathcal{C}$ the category of ($L_0^I$-complete) $H^*\widetilde{BN}$-modules. □

**Lemma 9.5.** *(Eilenberg-Moore) Consider the family $[\subseteq N] = \{H \leq G \mid H \subseteq N\}$ and the Quillen adjunction*
$$(-)^N : \mathrm{Mod}_{DEG_+}(\mathrm{Sp}_G) \leftrightarrows \mathrm{Mod}_{D\widetilde{BN}_+}(\mathrm{Sp}_W) : DEG_+ \otimes_{D\widetilde{BN}_+} -$$
*where we set $D\widetilde{BN}_+ = (DEG_+)^N$. Then for all $DEG_+$-modules $Y$, the counit map*
$$\epsilon_Y : DEG_+ \otimes_{D\widetilde{BN}_+} Y^N \to Y$$
*is a $E[\subseteq N]_+$-equivalence.*

*Proof.* A map of $G$-spectra is an $E[\subseteq N]_+$-equivalence if and only if its restriction to $N$-spectra is a weak equivalence. Therefore, it is sufficient to check that $DEN_+ \otimes_{DBN_+} Y^N \to Y$ is a weak equivalence. The full subcategory of $DEN_+$-module spectra $Y$ for which $\epsilon_Y$ is a weak equivalence is localizing and clearly contains $DEN_+$. Since $DEN_+$ generates $\mathrm{Mod}_{DEN_+}$ by [20, 3.1] the claim follows. □

We now ready to prove our main result.

**Theorem 9.6.** *Let $G$ be a compact Lie group with identity component $N$ and component group $W = G/N$. Let $I$ be the augmentation ideal of $H^*BN$. Then there is a symmetric monoidal Quillen equivalence*
$$L_{EG_+}(\mathrm{Sp}_G) \simeq_Q \mathrm{Mod}^{\wedge}_{H^*\widetilde{BN}[W]}$$
*between rational cofree $G$-spectra and $L_0^I$-complete dg-$H^*\widetilde{BN}[W]$-modules. In particular, there is a tensor-triangulated equivalence*
$$\textit{cofree G-spectra} \simeq_{\triangle} \mathcal{D}(L_0^I\textit{-complete } H^*\widetilde{BN}[W]\textit{-modules}).$$

*Proof.* We will prove the theorem using the Compactly Generated Localization Principle 3.14. To have a better control on the compact generators of the localized categories, it is convenient to change our model for cofree $G$-spectra. Thus we note that
$$L_{EG_+}\mathrm{Sp}_G = L_{G_+}\mathrm{Sp}_G$$
since the $EG_+$-equivalences are the same as the $G_+$-equivalences. Using Proposition 4.4 we have a symmetric monoidal Quillen equivalence $L_{G_+}(\mathrm{Sp}_G) \simeq_Q L_{G_+}(\mathrm{Mod}_{DEG_+})$.

Taking categorical $G$-fixed points loses too much information since $\mathrm{Mod}_{DEG_+}$ is no longer generated by $DEG_+$. Instead we slightly modify the model structure and then take $N$-fixed points. Consider the family $[\subseteq N] = \{H \leq G \mid H \subseteq N\}$. There is a symmetric monoidal Quillen equivalence
$$L_{G_+}\mathrm{Mod}_{DEG_+} \rightleftarrows L_{G_+}L_{E[\subseteq N]_+}\mathrm{Mod}_{DEG_+}$$
since $G_+ \wedge E[\subseteq N]_+ \to G_+$ is a weak equivalence.

We now take categorical $N$-fixed points to remove equivariance. We use the tilde in $D\widetilde{BN}_+ = (DEG_+)^N$ to emphasize that it may have a non-trivial $W$-action. We apply the Compactly Generated Localization Principle to the symmetric monoidal Quillen adjunction
$$(-)^N : L_{E[\subseteq N]_+}\mathrm{Mod}_{DEG_+}(\mathrm{Sp}_G) \leftrightarrows \mathrm{Mod}_{D\widetilde{BN}_+}(\mathrm{Sp}_W) : DEG_+ \otimes_{D\widetilde{BN}_+} -$$
to obtain a symmetric monoidal Quillen equivalence after localization. There are several conditions that need to be checked. Firstly, we claim that $L_{G_+}L_{E[\subseteq N]_+}\mathrm{Mod}_{DEG_+}$ is compactly generated by $DG_+ \simeq DG_+ \wedge DEG_+$. It is clear that it generates so we only show that it is compact. By definition of sum in the localized category, we have to show that

(1) $$\mathrm{hMod}_{DEG_+}(DG_+, F(EG_+, \bigvee_i Y_i)) \simeq \bigoplus_i \mathrm{hMod}_{DEG_+}(DG_+, Y_i)$$

where $Y_i$ is cofree for all $i$. This is now clear since $DG_+$ is small and $DG_+ \wedge EG_+ \simeq DG_+$. We also claim that $(DG_+)^N \simeq W_+$ compactly generates $L_{W_+}\mathrm{Mod}_{D\widetilde{BN}_+}(\mathrm{Sp}_W)$. Since $W_+$ has a trivial $D\widetilde{BN}_+$-action, it builds $D\widetilde{BN}_+ \wedge W_+$ in $\mathrm{Mod}_{D\widetilde{BN}_+}$ and hence it generates $L_{W_+}\mathrm{Mod}_{D\widetilde{BN}_+}(\mathrm{Sp}_W)$. It is also compact by a similar argument to (1). By the Compactly Generated Localization Principle it remains to check that the derived counit is a $G_+$-equivalence on $DG_+$, and that the derived counit is an $E[\subseteq N]_+$-equivalence for $G_+$. These are true by the Eilenberg-Moore Lemma. Hence we have a symmetric monoidal Quillen equivalence

$$L_{G_+} L_{E[\subseteq N]_+} \mathrm{Mod}_{DEG_+} \simeq_Q L_{W_+} \mathrm{Mod}_{D\widetilde{BN}_+}(\mathrm{Sp}_W).$$

Note that have an equality of model categories

$$L_{W_+} \mathrm{Mod}_{D\widetilde{BN}_+}(\mathrm{Sp}_W) = L_{W_+} \mathrm{Mod}_{D\widetilde{BN}_+}(L_{EW_+} \mathrm{Sp}_W)$$

since $EW_+ \wedge W_+ \simeq W_+$.

We can rewrite the target category as $L_{W_+} \mathrm{Mod}_{D\widetilde{BN}_+}(\mathrm{Sp}[W])$ and apply Shipley's theorem [45] to obtain symmetric monoidal Quillen equivalences

$$L_{W_+} \mathrm{Mod}_{D\widetilde{BN}_+}(\mathrm{Sp}[W]) \simeq_Q L_{\Theta(W_+)} \mathrm{Mod}_{\Theta D\widetilde{BN}_+}(\mathrm{Mod}_{\mathbb{Q}[W]}) \simeq_Q L_{\Theta(W_+)} \mathrm{Mod}_{\Theta D\widetilde{BN}_+[W]}.$$

One can construct a $\mathbb{Q}[W]$-module map $H^*\widetilde{BN} \to \Theta D\widetilde{BN}_+$ which is a quasi-isomorphism as in [26, §7]. Since the map is compatible with the $W$-action, there is a symmetric monoidal Quillen equivalence

$$\mathrm{Mod}_{\Theta D\widetilde{BN}_+[W]} \simeq_Q \mathrm{Mod}_{H^*\widetilde{BN}[W]}.$$

Note that $H_*(\Theta(W_+)) = H_0(\Theta(W_+)) = \mathbb{Q}[W]$ and hence $\Theta(W_+)$ is formal as a $H^*\widetilde{BN}[W]$-module since we have a zig-zag of quasi-isomorphisms $\Theta(W_+) \leftarrow \tau_{\geq 0}(\Theta(W_+)) \to H_0(\Theta(W_+))$ where $\tau_{\geq 0}$ denotes the connective cover functor. Putting all this together, we deduce a zig-zag of symmetric monoidal Quillen equivalences

$$L_{EG_+}(\mathrm{Sp}_G) \simeq_Q L_{\mathbb{Q}[W]} \mathrm{Mod}_{H^*\widetilde{BN}[W]}.$$

We now claim that $L_{\mathbb{Q}[W]} \mathrm{Mod}_{H^*\widetilde{BN}[W]} = (L_{\mathbb{Q}} \mathrm{Mod}_{H^*BN})[\widetilde{W}]$. As the underlying categories are equal and the acyclic fibrations are easily seen to be the same, we only need to argue that the model categories have the same weak equivalences. This is clear since

$$\mathbb{Q}[W] \otimes_{H^*\widetilde{BN}[W]} M \cong \mathbb{Q} \otimes_{H^*BN} M$$

for all $H^*\widetilde{BN}[W]$-modules $M$. Hence the two model categories are equal.

Finally, using Lemma 9.2 and Theorem 6.10, we conclude that there are symmetric monoidal Quillen equivalences

$$(L_{\mathbb{Q}} \mathrm{Mod}_{H^*BN})[\widetilde{W}] \simeq_Q \mathrm{Mod}^\wedge_{H^*BN}[\widetilde{W}] \simeq_Q \mathrm{Mod}^\wedge_{H^*\widetilde{BN}[W]}.$$

□

**Remark 9.7.** Our proof bridges a gap in [26]. In the cited paper it is stated that there is a Quillen equivalence

$$(-)^N : \mathrm{Cell}_{G_+} \mathrm{Mod}_{DEG_+}(\mathrm{Sp}_G) \leftrightarrows \mathrm{Cell}_{W_+} \mathrm{Mod}_{D\widetilde{BN}_+}(\mathrm{Sp}_W) : DEG_+ \otimes_{D\widetilde{BN}_+} -$$

obtained by the Cellularization Principle. The claim as it is stated it is not correct. Indeed, if we want to apply the Cellularization Principle we need to check that the counit $DEG_+ \otimes_{D\widetilde{BN}_+} (G_+)^N \to G_+$ is a weak equivalence of $G$-spectra, which in general is false. Nonetheless, we can modify the argument as follows. Firstly there is a Quillen equivalence

$$\mathrm{Cell}_{G_+} \mathrm{Mod}_{DEG_+} \rightleftarrows \mathrm{Cell}_{G_+} L_{E[\subseteq N]_+} \mathrm{Mod}_{DEG_+}.$$

We note that the localization $\mathrm{Cell}_{G_+} L_{E[\subseteq N]_+} \mathrm{Mod}_{DEG_+}$ exists, since left Bousfield localizations of right proper, stable model categories are right proper by [10, 4.7]. We can then apply the Cellularization Principle to the Quillen adjunction

$$(-)^N : L_{E[\subseteq N]_+} \mathrm{Mod}_{DEG_+}(\mathrm{Sp}_G) \leftrightarrows \mathrm{Mod}_{D\widetilde{BN}_+}(\mathrm{Sp}_W) : DEG_+ \otimes_{D\widetilde{BN}_+} -$$

and the Eilenberg-Moore Lemma to show that this is a Quillen equivalence after cellularization.

## 10. Adams spectral sequence

In this section, we construct an Adams spectral sequence for cofree $G$-spectra. It provides a tool for calculating the space of maps between two cofree $G$-spectra in terms of $L_0^I$-complete modules, and furthermore gives intuition for the Quillen equivalence given in the previous section.

We describe the construction of an Adams spectral sequence based on projective resolutions as in [2].

Let $\mathcal{T}$ be a triangulated category and let $\mathcal{A}$ be a $\mathbb{Z}$-graded abelian category with enough projectives. Note that $\mathcal{T}(X,Y)$ is a $\mathbb{Z}$-graded abelian group via $\mathcal{T}(X,Y)_n = \mathcal{T}(\Sigma^n X, Y)$. Assume that we are given a $\mathbb{Z}$-graded exact functor $\pi_*^{\mathcal{A}} \colon \mathcal{T} \to \mathcal{A}$. We aim to construct a conditionally convergent Adams-type spectral sequence
$$E_2^{s,t} = \operatorname{Ext}_{\mathcal{A}}^{s,t}(\pi_*^{\mathcal{A}}(X), \pi_*^{\mathcal{A}}(Y)) \Rightarrow \mathcal{T}(X,Y)_{t-s}$$
for all $X, Y \in \mathcal{T}$. We list the steps needed.

**Step 0:** Choose a projective resolution of $\pi_*^{\mathcal{A}}(X)$ in $\mathcal{A}$
$$0 \leftarrow \pi_*^{\mathcal{A}}(X) \leftarrow P_0 \leftarrow P_1 \leftarrow P_2 \leftarrow \ldots.$$

**Step 1:** Realize the projectives, i.e., find $\mathbb{P}_j \in \mathcal{T}$ so that $\pi_*^{\mathcal{A}}(\mathbb{P}_j) = P_j$.

**Step 2:** Let $X \in \mathcal{T}$ and $\mathbb{P}_j$ as above. Show that the functor $\pi_*^{\mathcal{A}}$ induces an isomorphism
$$\mathcal{T}(\mathbb{P}_j, X) \xrightarrow{\cong} \operatorname{Hom}_{\mathcal{A}}(P_j, \pi_*^{\mathcal{A}}(X)).$$

**Step 3:** Using Step 0 and Step 1, we can formally produce a tower

$$\begin{array}{ccccccc}
\mathbb{P}_0 & & \Sigma^1 \mathbb{P}_1 & & \Sigma^2 \mathbb{P}_2 & & \Sigma^3 \mathbb{P}_3 \\
\downarrow & & \downarrow & & \downarrow & & \downarrow \\
X = X_0 & \longrightarrow & X_1 & \longrightarrow & X_2 & \longrightarrow & X_3 & \longrightarrow \cdots
\end{array}$$

**Step 4:** Apply $\mathcal{T}(-, Y)$ to get a spectral sequence with $E_1$-page:
$$E_1^{s,*} = \mathcal{T}(\mathbb{P}_s, Y) = \operatorname{Hom}_{\mathcal{A}}(P_s, \pi_*^{\mathcal{A}}(Y)).$$

By construction, we will have a conditionally convergent spectral sequence
$$E_2^{*,*} = \operatorname{Ext}_{\mathcal{A}}^{*,*}(\pi_*^{\mathcal{A}}(X), \pi_*^{\mathcal{A}}(Y)) \Rightarrow \mathcal{T}(\overline{X}, Y)_*$$
where $\overline{X}$ is the fibre of the canonical map $X \to \operatorname{hocolim}_s X_s$.

**Step 5:** Show that $\operatorname{hocolim}_s X_s = 0$.

We apply the recipe above in the following setting. Fix a compact Lie group $G$ with identity component $N$ and component group $W$, and fix $I$ to be the augmentation ideal of $H^*BN$. We consider the abelian category $\operatorname{Mod}_{H^*\widetilde{BN}[W]}^{\wedge}$ of graded $L_0^I$-complete modules over the skewed group ring $H^*\widetilde{BN}[W]$, and the homotopy category of rational cofree $G$-spectra. Before we give the exact functor, we need a preliminary result.

**Remark 10.1.** We recall a spectral sequence relating local homology to equivariant homotopy groups, see [23]. Let $R$ be a ring $G$-spectrum and $M$ an $R$-module. For $J = (x_1, \ldots, x_r)$ a finitely generated ideal in $\pi_*^G R$ define
$$M_J^{\wedge} = F(K(J), M)$$
where $K(J) = \operatorname{fib}(R \to R[1/x_1]) \otimes_R \cdots \otimes_R \operatorname{fib}(R \to R[1/x_r])$ is the Koszul spectrum. Then there is a convergent spectral sequence
$$E_2^{*,*} = L_*^J(\pi_*^G R; \pi_*^G M) \implies \pi_*^G(M_J^{\wedge}).$$
In the special case that $R$ has Thom isomorphisms and $J$ is the augmentation ideal of $\pi_*^G R$, there is an equivalence $M_J^{\wedge} \xrightarrow{\sim} F(EG_+, M)$ by [23, 2.5].

**Lemma 10.2.** *Let $X$ be a cofree $G$-spectrum. Then $\pi_*^N X$ is $L_0^I$-complete.*

*Proof.* For any $DEN_+$-module $M$, there is a convergent spectral sequence
$$E_2^{*,*} = L_*^I(H^*BN; \pi_*^N M) \implies \pi_*^N(F(EG_+, M))$$
by Remark 10.1. If in addition $M$ is cofree, we also have that $M_I^\wedge \simeq M$.

Now let $X$ be a cofree $G$-spectrum. The discussion above tells us that we have a convergent spectral sequence
$$E_2^{*,*} = L_*^I(H^*BN; \pi_*^N X) \implies \pi_*^N X.$$

Since the $E_2$-page of the spectral sequence consists of $L_0^I$-complete modules by [21, 4.1], and the kernel and cokernel of a map of $L_0^I$-complete modules is $L_0^I$-complete, we have that $\pi_*^N X$ is $L_0^I$-complete.

Finally, note that $W$ acts on $\pi_*^N(X)$ by conjugation, making it naturally a module over $H^*\widetilde{BN}[W]$. $\square$

Therefore we may use the exact functor
$$\pi_*^N \colon \mathrm{hSp}_G^{\mathrm{cofree}} \to \mathrm{Mod}^\wedge_{H^*\widetilde{BN}[W]}$$
for the construction of the Adams spectral sequence.

**Lemma 10.3** (Step 1). *The abelian category $\mathrm{Mod}^\wedge_{H^*\widetilde{BN}[W]}$ has enough projectives. Moreover, the projectives are realized, that is*
$$\pi_*^N(F(EG_+, \bigvee \Sigma^{n_i} S^0) \wedge W_+) \cong L_0^I(\bigoplus \Sigma^{n_i} H^*\widetilde{BN})[W].$$

*Proof.* Using that $L_0^I$ is right exact and left adjoint to the inclusion, we see that $L_0^I(\bigoplus \Sigma^{n_i} H^*\widetilde{BN})[W]$ is projective in $\mathrm{Mod}^\wedge_{H^*\widetilde{BN}[W]}$ and that there are enough projectives. It is left to show that the projectives are realized. Note that
$$\pi_*^N(F(EG_+, \bigvee \Sigma^{n_i} S^0) \wedge W_+) \cong \pi_*^N(F(EG_+, \bigvee \Sigma^{n_i} S^0))[W]$$
so it is enough to show that
$$\pi_*^N(F(EN_+, \bigvee \Sigma^{n_i} S^0) \cong L_0^I(\oplus \Sigma^{n_i} H^*BN).$$

By isotropy separation,
$$F(EN_+, \bigvee \Sigma^{n_i} S^0) \xrightarrow{\sim} F(EN_+, \bigvee \Sigma^{n_i} DEN_+).$$

There is a spectral sequence
$$E_2^{*,*} = L_*^I(\pi_*^N M) \implies [EN_+, M]_*^N = \pi_*^N(F(EN_+, M))$$
and when $M = \bigvee \Sigma^{n_i} DEN_+$ the $E_2$-page has the form
$$L_*^I \pi_*^N(\bigvee \Sigma^{n_i} DEN_+) \cong L_*^I(\oplus \Sigma^{n_i} H^*BN) = L_0^I(\oplus \Sigma^{n_i} H^*BN)$$
by Proposition 6.2. Since the $E_2$-page is concentrated in one line, the spectral sequence collapses so that
$$L_0^I(\oplus \Sigma^{n_i} H^*BN) \cong \pi_*^N(F(EN_+, \bigvee \Sigma^{n_i} DEN_+)).$$
$\square$

We also need to realize the maps.

**Lemma 10.4** (Step 2). *Taking homotopy groups gives an isomorphism*
$$\pi_*^N \colon [F(EG_+, \bigvee \Sigma^{n_i} S^0) \wedge W_+, Y]_*^G \xrightarrow{\cong} \mathrm{Hom}_{H^*\widetilde{BN}[W]}(L_0^I(\oplus \Sigma^{n_i} H^*\widetilde{BN})[W], \pi_*^N(Y))^4.$$

---

[4] The proof given here relies upon Corollary 5.7 which is incorrect. However the lemma is true as stated; see the corrigendum at the end of this paper for an alternative proof

*Proof.* We apply the change of groups adjunctions on both sides to reduce to showing that
$$\pi_*^N \colon [F(EN_+, \bigvee \Sigma^{n_i} S^0), Y]_*^N \xrightarrow{\cong} \mathrm{Hom}_{H^*BN}(\oplus \Sigma^{n_i} H^*BN, \pi_*^N(Y)).$$
Since there is a weak equivalence $EN_+ \wedge X \xrightarrow{\sim} EN_+ \wedge F(EN_+, X)$ for any $X$, we see that
$$[F(EN_+, \bigvee \Sigma^{n_i} S^0), F(EN_+, Y)]^N \cong [\bigvee F(EN_+, \Sigma^{n_i} S^0), F(EN_+, Y)]^N.$$
Accordingly, it is enough to show that
$$\pi_*^N \colon [DEN_+, Y]_*^N \xrightarrow{\cong} \mathrm{Hom}_{H^*BN}(H^*BN, \pi_*^N(Y))$$
for all $Y$ cofree $N$-spectra. We show that both sides are homology theories on cofree $N$-spectra, taking values in $L_0^I$-complete $H^*BN$-modules. In other words, we show that they are exact and satisfy the wedge axiom. Exactness follows from the fact that $\pi_*^N$ sends triangles of cofree $N$-spectra to long exact sequences of $L_0^I$-complete modules by Lemma 10.2. We now turn to the wedge axiom. By definition of coproducts in the categories of cofree spectra and $L_0^I$-complete modules, we need to show that
$$L_0^I(\oplus \pi_*^N X_i) \cong \pi_*^N(F(EG_+, \bigvee X_i))$$
for $X_i$ cofree $N$-spectra. By Remark 10.1 there is a convergent spectral sequence
$$L_*^I(\oplus \pi_*^N X_i) \implies \pi_*^N(F(EG_+, \bigvee X_i)).$$
This spectral sequence collapses by Corollary 5.7 and Lemma 10.2, which proves the wedge axiom.

As both the left and right hand sides are homology theories on the homotopy category of cofree $N$-spectra, it is enough to show that the natural transformation $\pi_*^N$ induces an isomorphism for the compact generator $Y = DN_+$. It is easy to see that both homology theories evaluated at $DN_+$ give $\mathbb{Q}$. Thus we only have to argue that the natural transformation $\pi_*^N$ is non-zero. Observe that a nontrivial $N$-equivariant map $f \colon DEN_+ \to DN_+$ corresponds to a nontrivial $N$-equivariant map $\widetilde{f} \colon N_+ \to EN_+$ which gives a map $\widetilde{f}/N_+ \colon S^0 \to BN_+$ which is nontrivial in reduced $H_0$. It remains to note that for a free $N$-spectrum we have $\pi_*^N(X) = H_*(X/N)$ up to an integer shift. $\square$

Since $\pi_*^N$ maps homotopy direct limits to direct limits, it is left to show the following:

**Lemma 10.5** (Step 5). *Let $X$ be a cofree $G$-spectrum with $\pi_*^N(X) = 0$. Then $X \simeq 0$.*

*Proof.* We first prove the claim for the connected case and then we show how to extend it to all compact Lie groups. Note that there is an equivalence $EN_+ \simeq EN_+ \wedge DEN_+$ so that $EN_+ \in \mathrm{Mod}_{DEN_+}$. We claim that $N_+ \in \mathrm{Loc}(EN_+)$; that is $N_+$ is in the localizing subcategory generated by $EN_+$. Note that $DN_+$ is cofree and hence a $DEN_+$-module. Since $DEN_+$ generates the category $\mathrm{Mod}_{DEN_+}$ by [20, 3.1], we get that $DN_+ \in \mathrm{Loc}_{\mathrm{Mod}_{DEN_+}}(DEN_+)$. Since the forgetful functor $\mathrm{Mod}_{DEN_+} \to \mathrm{Sp}_N$ and $EN_+ \wedge - \colon \mathrm{Sp}_N \to \mathrm{Sp}_N$ preserve colimits we get $EN_+ \wedge DN_+ \in \mathrm{Loc}(EN_+ \wedge DEN_+)$. By the Wirthmüller Isomorphism, we see that $DN_+ \simeq \Sigma^{-d} N_+$ where $d$ is the dimension of $N$. Putting all this together, $N_+ \in \mathrm{Loc}(EN_+)$ as required. Let us now prove that for a cofree $N$-spectrum $X$ with $\pi_*^N(X) = 0$, then $\pi_*(X) = 0$ and hence $X \simeq 0$. By hypothesis, we have
$$0 = \pi_*^N(X) = [EN_+, X]^N.$$
By a localizing subcategory argument we get $\pi_*(X) = [N_+, X]^N = 0$ as required. Finally, let $G$ be any compact Lie group and let $X$ be a cofree $G$-spectrum with $\pi_*^N(X) = 0$. By the previous paragraph, we know that $X$ is $N$-equivariantly contractible, that is $F(W_+, X) \simeq 0$ and hence $F(EW_+, X) \simeq 0$. Therefore
$$X \simeq F(\widetilde{E}W, X) \simeq F(\widetilde{E}W \wedge EG_+, X) \simeq 0$$
since $X$ is cofree. $\square$

Finally, we have our Adams spectral sequence:

**Theorem 10.6.** *For $X$ and $Y$ cofree $G$-spectra, there is a strongly convergent Adams spectral sequence*
$$E_2^{*,*} = \mathrm{Ext}^{*,*}_{H^*\widetilde{BN}[W]}(\pi_*^N X, \pi_*^N Y) \implies [X, Y]_*^G.$$

*Proof.* Combining the results of this section with Proposition 5.6, we have constructed a conditionally convergent spectral sequence as above. Note that $H^*\widetilde{BN}[W]$ has global homological dimension smaller or equal to $r = \mathrm{rank}(N)$ since the projectives are induced from the category of $H^*BN$-modules, see Lemma 10.3. It follows that the spectral sequence is concentrated in rows 0 to $r$, and hence is strongly convergent. □

(Pol) SCHOOL OF MATHEMATICS AND STATISTICS, HICKS BUILDING, SHEFFIELD S3 7RH, UK

*Email address*: `lpol1@sheffield.ac.uk`

(Williamson) SCHOOL OF MATHEMATICS AND STATISTICS, HICKS BUILDING, SHEFFIELD S3 7RH, UK

*Email address*: `jwilliamson3@sheffield.ac.uk`

# CORRIGENDUM TO "THE LEFT LOCALIZATION PRINCIPLE, COMPLETIONS, AND COFREE $G$-SPECTRA"

LUCA POL AND JORDAN WILLIAMSON

There are two errors in the paper [3]. The errors are in certain auxiliary arguments, and we show here that all of the main results of our paper remain correct as stated. Specifically, the only numbered results of [3] that need change are 3.4 (2), 3.7 (2) and 5.7. In particular, Lemma 10.4 remains true as stated, and we give a new proof here. We are grateful to Neil Strickland for drawing these errors to our attention, and to John Greenlees for various discussions about these.

The first error occurs in Theorem 3.4(2) (and is repeated in Theorem 3.7(2)). Under the stated hypotheses, the derived functor $GR$ need not preserve sums and so one cannot proceed by a localizing subcategory argument as suggested in the proof. To fix this, the assumption that the derived counit $FQGRL \to RL$ is a $\mathsf{T}$-equivalence for all $L \in \mathcal{L}$ needs to be strengthened to the assumption that the derived counit is a weak equivalence in $\mathcal{D}$ for all $L \in \mathcal{L}$. Under this stronger assumption $GR$ preserves sums and the proof proceeds as in the proof for Part (1). Alternatively, under this stronger assumption Part (2) is a consequence of Part (1) by setting $S = GRT$. We note that this does not affect any other results in the paper.

Before discussing the second mistake, we point out a possible source of confusion. Let $R$ be a graded commutative ring and $I$ a finitely generated homogeneous ideal of $R$. For a graded $R$-module $M$, we write $L_n^I M$ for the $n$-th left derived functor of $I$-adic completion applied to $M$. By applying $L_n^I$ to the underlying graded module, we can extend this to an endofunctor of dg-modules. Note that this is in line with the definition of $L_0^I$-complete dg-modules given in Definition 5.1. In particular, this functor should *not* be confused with the $n$-th homology of the total left derived functor of $I$-adic completion, which is instead computed using dg-projective resolutions. The only appearance of the latter interpretation is at the end of Theorem 5.5, and plays no role in the rest of the paper.

The second error is Corollary 5.7 which is not true in the generality stated. A counterexample is given by taking $R = k[x,y]$, $I = (x,y)$, $M_i = R/(a^i, b^i)$ and $P_i$ to be the Koszul resolution of $M_i$. If the ideal $I$ is principal, then Corollary 5.7 remains true since in this case $L_1^I$ is left exact by [1, 1.9] and so $L_1^I(\bigoplus M_i) \to L_1^I(\prod M_i) = 0$ is a monomorphism. More generally, if $I$ is generated by a sequence of length $n$, then $L_m^I(\bigoplus M_i)$ is zero for all $m \geq n$, and the above counterexample shows that there always exists a collection $M_i$ for which $L_{n-1}^I(\bigoplus M_i)$ is not zero.

The error in the statement of Corollary 5.7 only impacts the proof of Lemma 10.4. Here we give an alternative proof of Lemma 10.4 which does not rely upon Corollary 5.7. As it may be of independent interest we state the new key ingredient separately. Firstly, we need to recall some terminology.

Let $\mathsf{T}$ be a triangulated category and let $X$ be an object of $\mathsf{T}$; for us, $\mathsf{T}$ will be the (homotopy) category of rational cofree $G$-spectra. A full replete subcategory of $\mathsf{T}$ which is closed under retracts is said to be thick. If in addition it is closed under arbitrary coproducts (resp. products) it is said to be localizing (resp. colocalizing). The thick (resp. localizing, resp. colocalizing) subcategory of $\mathsf{T}$ generated by $X$ is the smallest thick (resp. localizing, resp. colocalizing) subcategory of $\mathsf{T}$ which contains $X$. We then say that an object $Y$ of $\mathsf{T}$ is:

- *finitely built* from $X$ if $Y$ is in the thick subcategory generated by $X$;
- *built* from $X$ if $Y$ is in the localizing subcategory generated by $X$;
- *cobuilt* from $X$ if $Y$ is in the colocalizing subcategory generated by $X$;

**Lemma 1.** *Let $G$ be a connected compact Lie group.*

(a) *Rationally, $EG_+$ finitely builds $G_+$.*
(b) *The rational spectrum $\bigoplus_i \Sigma^{n_i} G_+$ is a retract of $\prod_i \Sigma^{n_i} G_+$.*
(c) *The category of rational cofree $G$-spectra coincides with the colocalizing subcategory generated by $DEG_+$, and with the colocalizing subcategory generated by $DG_+$.*

*Proof.* The proof of (a) can be found in [2, 6.4]. Here we repeat the argument for the benefit of the reader. Since $G$ is connected, we know that $H^*(BG; \mathbb{Q}) = \mathbb{Q}[x_1, \ldots, x_r]$ with generators in even degrees, where $r = \mathrm{rank}(G)$. Note that $[EG_+, EG_+]_*^G = H^*(BG; \mathbb{Q})$ so we can view each generator $x_i$ as an endomorphism of $EG_+$. By the Wirthmüller isomorphism, we have that $\pi_*^G(G_+) = \Sigma^d \mathbb{Q}$ where $d$ is the dimension of $G$. This suggests that we can build $G_+$ from $EG_+$ via a Koszul resolution. Define

$$EG_+//x_1 = \mathrm{fib}(\Sigma^{-|x_1|} EG_+ \xrightarrow{x_1} EG_+)$$

and

$$EG_+//(x_1, \ldots x_j) = \mathrm{fib}(\Sigma^{-|x_j|} EG_+//(x_1, \ldots, x_{j-1}) \xrightarrow{x_j} EG_+//(x_1, \ldots, x_{j-1}))$$

for $1 < j \leq r$. By construction, $EG_+$ finitely builds $EG_+//(x_1, \ldots, x_r)$. The latter spectrum can be identified with the bottom cell of $EG_+$ which is $G_+$. This concludes the proof of (a).

For part (b), we first note that in the category of $\mathbb{Q}$-vector spaces the direct sum is a retract of the product. It follows that the rational (non-equivariant) spectrum $\bigoplus S^{n_i}$ is a retract of $\prod S^{n_i}$. The claim in (b) then follows by applying the functor $G_+ \wedge -$ to this splitting together with the identification

$$G_+ \wedge (\prod S^{n_i}) \simeq \Sigma^d F(G_+, \prod S^{n_i}) \simeq \prod \Sigma^{n_i} G_+$$

using that $G_+ \simeq \Sigma^d DG_+$.

We now prove part (c). Note that $G_+$ builds $EG_+$ as it is free. It follows that $DG_+$ cobuilds $DEG_+$, and so to prove part (c) it is sufficient to show that $DEG_+$ cobuilds any rational cofree $G$-spectrum $Y$. We have seen in [3, 10.3] that any free resolution of $\pi_*^G(Y)$ can be realized in rational cofree $G$-spectra by $\widehat{\mathbb{F}}_i = F(EG_+, \bigvee \Sigma^{n_i} DEG_+)$. Since the homological dimension of the category of $L_0^I$-complete $H^*BG$-modules is $r = \mathrm{rank}(G)$, it follows that $Y$ is finitely built from the collection of $\widehat{\mathbb{F}}_i$. Therefore it suffices to show that each $\widehat{\mathbb{F}}_i$ is cobuilt from $DEG_+$.

We write $\mathbb{F} = \bigvee \Sigma^{n_i} DEG_+$ and $\widehat{\mathbb{F}} = F(EG_+, \bigvee \Sigma^{n_i} DEG_+)$ for its cofree-ification. Since $EG_+$ finitely builds $G_+$, applying the dual functor shows that $DEG_+$ cobuilds $DG_+ \simeq \Sigma^{-d} G_+$ (by the Wirthmüller isomorphism). Note that $DEG_+ \wedge G_+ \simeq G_+$, and therefore $DEG_+$ cobuilds $\prod(\Sigma^{n_i} DEG_+ \wedge G_+)$. Hence $DEG_+$ cobuilds $\mathbb{F} \wedge G_+ = \bigvee \Sigma^{n_i} DEG_+ \wedge G_+$ since it is a retract of the product by part (b).

As $EG_+^{(n)}$ is a finite free $G$-CW-complex, $\mathbb{F} \wedge G_+$ finitely builds $\mathbb{F} \wedge DEG_+^{(n)}$. Therefore $DEG_+$ cobuilds $\mathbb{F} \wedge DEG_+^{(n)}$. Finally, we note that $\widehat{\mathbb{F}} = \mathrm{holim}(\mathbb{F} \wedge DEG_+^{(n)})$ since $EG_+^{(n)}$ is finite. This completes the proof of part (c). □

We restate Lemma 10.4 for completeness. Recall that for $G$ a compact Lie group, we write $N$ for the identity component and $W = G/N$ for the component group.

**Lemma 2.** *Taking homotopy groups gives an isomorphism*

$$\pi_*^N \colon [F(EG_+, \bigvee \Sigma^{n_i} S^0) \wedge W_+, Y]_*^G \to \mathrm{Hom}_{H^*\widetilde{BN}[W]}(L_0^I(\oplus \Sigma^{n_i} H^* \widetilde{BN})[W], \pi_*^N(Y)).$$

*Proof.* As stated in the proof given in the paper [3, Proof of 10.4], it is enough to show that

$$\pi_*^N \colon [DEN_+, Y]_*^N \to \mathrm{Hom}_{H^*BN}(H^*BN, \pi_*^N(Y))$$

is an isomorphism for all cofree $N$-spectra $Y$. We note that the collection of $Y$ for which the above map is an isomorphism is colocalizing. In [3, Proof of 10.4] we argued that the map is an isomorphism for $Y = DN_+$. As every rational cofree $N$-spectrum is cobuilt from $DN_+$ by Lemma 1 the claim follows. □

(Pol) Fakultät für Mathematik, Universität Regensburg, Universitätsstrasse 31, 93053 Regensburg, Germany

*Email address*: `luca.pol@ur.de`

(Williamson) Department of Algebra, Faculty of Mathematics and Physics, Charles University in Prague, Sokolovská 83, 186 75 Praha, Czech Republic

*Email address*: `williamson@karlin.mff.cuni.cz`